\input epsf
\documentstyle{amsppt}
\pagewidth{5.4truein}\hcorrection{0.55in}
\pageheight{7.5truein}\vcorrection{0.75in}
\TagsOnRight
\NoRunningHeads
\catcode`\@=11
\def\logo@{}
\footline={\ifnum\pageno>1 \hfil\folio\hfil\else\hfil\fi}
\topmatter
\title Another dual of MacMahon's theorem on plane partitions
\endtitle
\author Mihai Ciucu\endauthor
\thanks Research supported in part by NSF grants DMS-1101670 and the Galileo Galilei Institute of Physics, Firenze, Italy.
\endthanks
\affil
  Department of Mathematics, Indiana University\\
  Bloomington, Indiana 47405
\endaffil
\abstract
In this paper we introduce a counterpart structure to the shamrocks studied in the paper {\it A dual of Macmahon's theorem on plane partitions} by M. Ciucu and C. Krattenthaler ({\it Proc. Natl. Acad. Sci. USA} {\bf 110} (2013), 4518--4523), which, just like the latter, can be included at the center of a lattice hexagon on the triangular lattice so that the region obtained from the hexagon by removing it has its number of lozenge tilings given by a simple product formula. The new structure, called a fern, consists of an arbitrary number of equilateral triangles of alternating orientations lined up along a lattice line. The shamrock and the fern seem to be the only structures with this property. It would be interesting to understand why these are the only two such structures.
\endabstract
\endtopmatter

\document

\def\mysec#1{\bigskip\centerline{\bf #1}\message{ * }\nopagebreak\bigskip\par}

\def\myref#1{\item"{[{\bf #1}]}"}

\def\epf{\hfill{$\square$}\smallpagebreak}

\def\cite#1{\relaxnext@
  \def\nextiii@##1,##2\end@{[{\bf##1},\,##2]}%
  \in@,{#1}\ifin@\def\next{\nextiii@#1\end@}\else
  \def\next{[{\bf#1}]}\fi\next}
\def\proclaimheadfont@{\smc}

\define\M{\operatorname{M}}

\define\h{\operatorname{H}}
%\define\M{\operatorname{M}}
%\define\M{\operatorname{M}}
%\define\M{\operatorname{M}}

\define\twoline#1#2{\line{\hfill{\smc #1}\hfill{\smc #2}\hfill}}
\define\twolinetwo#1#2{\line{{\smc #1}\hfill{\smc #2}}}
\define\twolinethree#1#2{\line{\phantom{poco}{\smc #1}\hfill{\smc #2}\phantom{poco}}}

\def\mypic#1{\epsffile{#1}}

%\def\epsfsize#1#2{0.36#1}
%\topinsert
%\centerline{\mypic{2-1.eps}}
%\centerline{{\smc Figure~2.1{\rm (a). $R_{(2,4,5),(2,4)}(2)$.}}}
%\endinsert

%\topinsert
%\twoline{\mypic{2-1b.eps}}{\mypic{2-1c.eps}}
%\twoline{Figure~2.1{\rm (b). $R_{\emptyset,(2,4)}(4)$.}}{Figure~2.1{\rm (c). 
%$R_{(2,4,5),\emptyset}(2)$.}}
%\endinsert

% ref nos
\define\And{1}
\define\cekz{2}
%\define\sc{3}
\define\ec{3}
\define\ov{4}
%\define\ef{6}
\define\vf{5}
\define\gd{6}
%\define\anglepap{8}
\define\CLP{7}
%\define\DT{10}
%\define\Glaish{12}
\define\GT{8}
\define\KKZ{9}
\define\Kuo{10}
\define\Kup{11}
\define\MacM{12}
\define\Sta{13}
\define\Ste{14}

% eq nos
\define\eaa{1.1}
\define\eab{1.2}
\define\eac{1.3}
\define\ead{1.4}
\define\eae{1.5}
\define\eaf{1.6}

\define\ebaa{2.1}
\define\ebab{2.2}
\define\eba{2.3}
\define\ebb{2.4}
\define\ebc{2.5}
\define\ebd{2.6}
\define\ebe{2.7}
\define\ebf{2.8}
\define\ebg{2.9}
\define\ebh{2.10}
\define\ebga{2.11}
\define\ebha{2.12}
\define\ebgb{2.13}
\define\ebhb{2.14}

\define\eca{3.1}
\define\ecb{3.2}
\define\ecc{3.3}
\define\ecd{3.4}
\define\ece{3.5}
\define\ecf{3.6}
\define\ecg{3.7}

\define\eci{3.9}
\define\ecj{3.10}
\define\eck{3.11}

\define\eda{4.1}
\define\edb{4.2}
\define\edc{4.3}
\define\edd{4.4}
\define\ede{4.5}
\define\edf{4.6}
\define\edg{4.7}
\define\edh{4.8}
\define\edi{4.9}
\define\edj{4.10}
\define\edk{4.11}
\define\edl{4.12}
\define\edm{4.13}

\define\eea{5.1}
\define\eeb{5.2}

% th nos
\define\taa{1.1}

\define\tba{2.1}

\define\tca{3.1}
\define\tcb{3.2}

% fig nos
\define\faa{1.1}
\define\fab{1.2}
\define\fac{1.3}
\define\fad{1.4}

\define\fba{2.1}
\define\fbb{2.2}

\define\fca{3.1}
\define\fcb{3.2}
\define\fcc{3.3}
\define\fcd{3.4}
\define\fce{3.5}

\define\fea{5.1}

%\vskip-0.1in
\mysec{1. Introduction}

Few results have influenced more the current intense attention devoted to different aspects of tilings and perfect matchings than MacMahon's classical theorem on the number of plane partitions that fit in a given box (see \cite{\MacM}, and \cite{\Sta}\cite{\And}\cite{\Kup}\cite{\Ste}\cite{\KKZ}\cite{\vf} for more recent developments). It is equivalent to the fact that the number of lozenge tilings of a hexagon of side-lengths $a$, $b$, $c$, $a$, $b$, $c$ (in cyclic order) on the triangular lattice is equal to
$$
P(a,b,c):=
\frac{\h(a)\h(b)\h(c)\h(a+b+c)}
{\h(a+b)\h(a+c)\h(b+c)},
\tag\eaa
$$
where the hyperfactorial $\h(n)$ is defined by
$$
\h(n):=0!\,1!\cdots(n-1)!
\tag\eab
$$
(see Figure {\faa} for an example). 

The striking elegance of this resut compels one to seek similar ones, which could help place it in a broader context. One step in this direction was taken in \cite{\vf}, where it was shown that if instead of the hexagon --- which is the region on the triangular lattice traced out by turning 60 degrees at each corner --- one considers the figure obtained by turning 120 degrees at each corner --- which we called a shamrock in \cite{\vf} (see Figure {\fab} for an illustration) --- then it is the {\it exterior} of that region which has, in a certain precise sense, a normalized number of tilings given by a simple product formula analogous to (\eaa). 

%\topinsert
%\twoline{\mypic{2-1b.eps}}{\mypic{2-1c.eps}}
%\twoline{Figure~2.1{\rm (b). $R_{\emptyset,(2,4)}(4)$.}}{Figure~2.1{\rm (c). 
%$R_{(2,4,5),\emptyset}(2)$.}}
%\endinsert

\topinsert
\twoline{\mypic{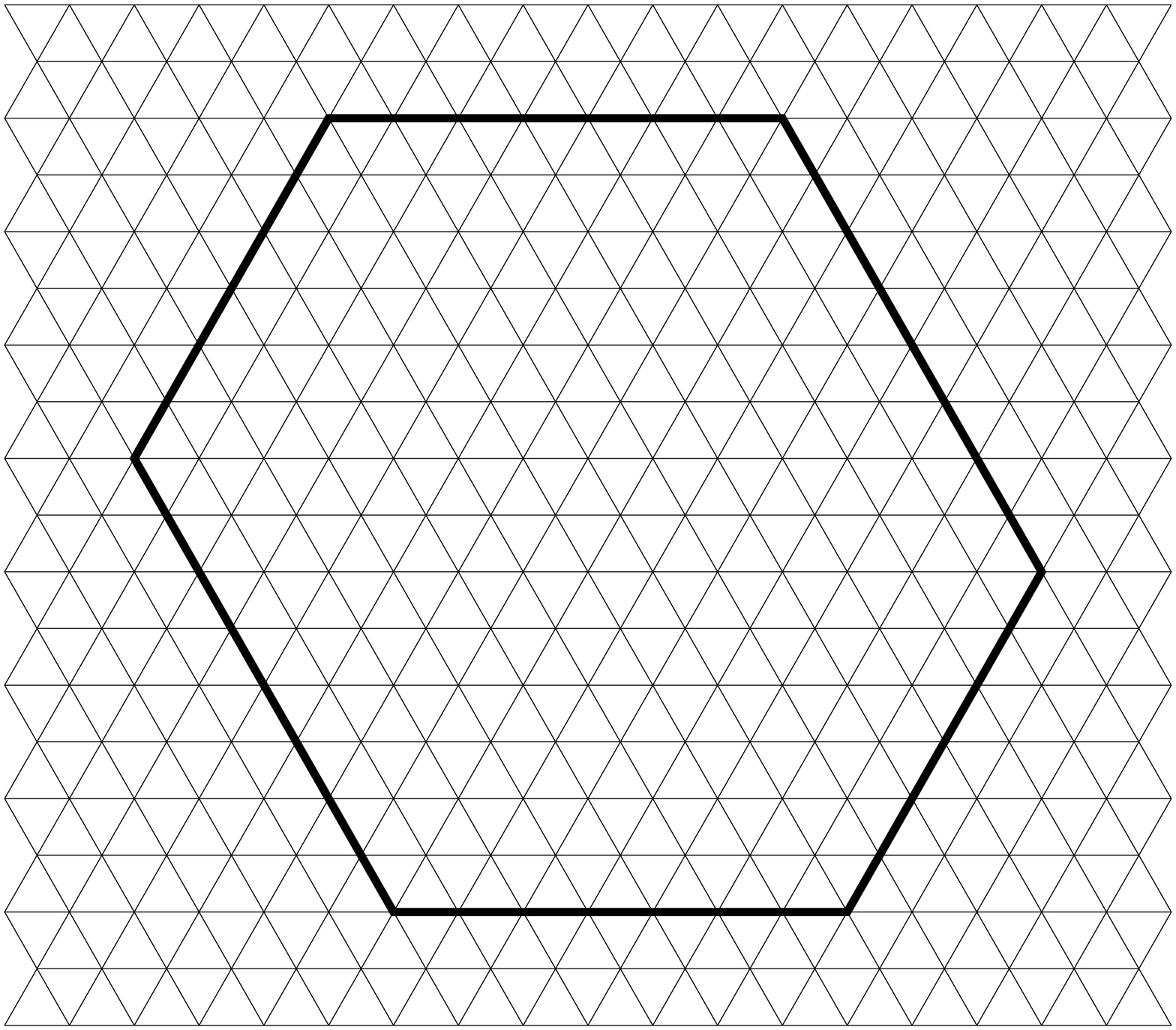}}{\mypic{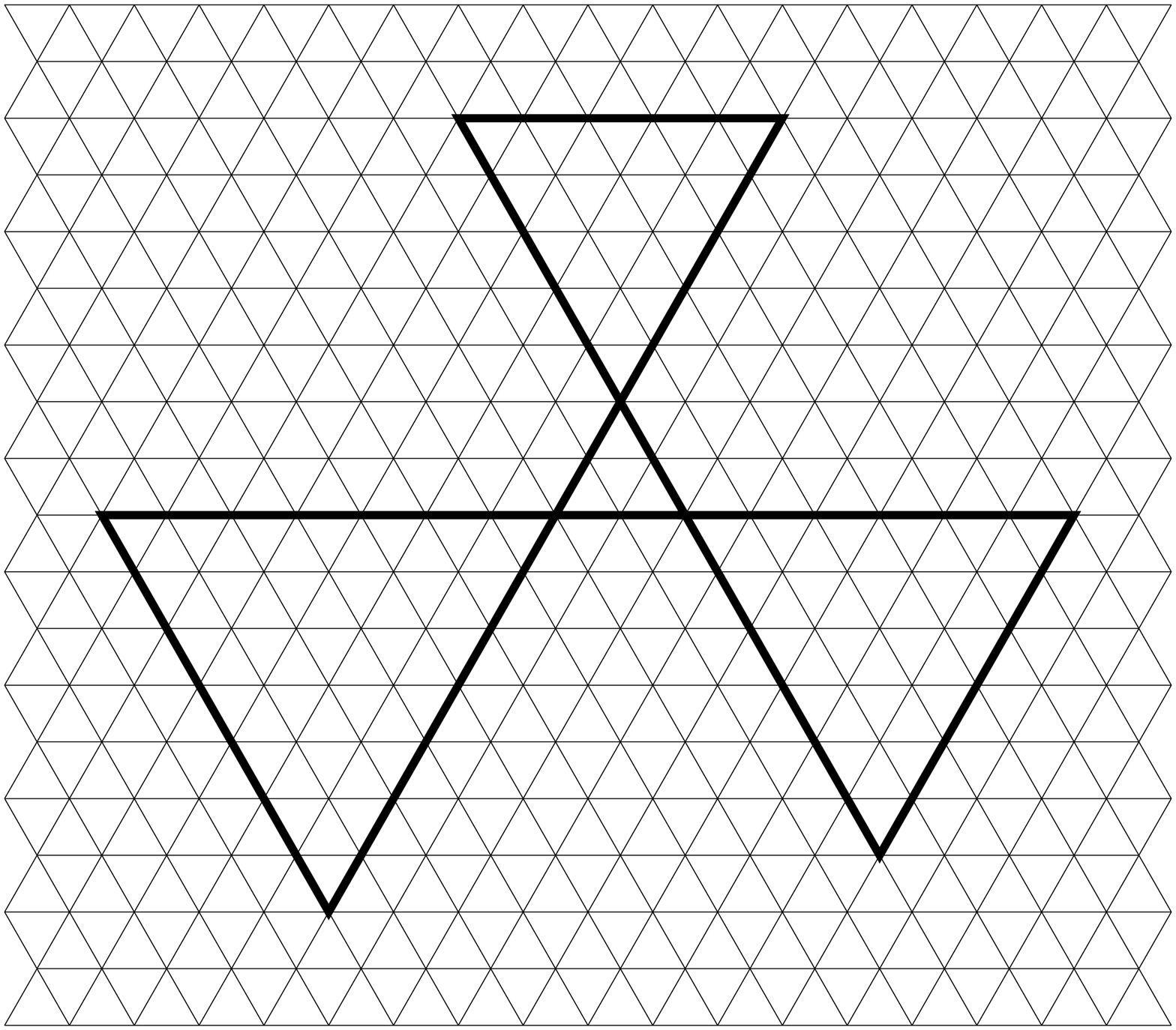}}
\twoline{Figure~{\faa}. {\rm Hexagon with $a=7$,}}{Figure~{\fab}. {\rm Shamrock with $m=2$,}}
\twoline{ {\rm \ \ \ \ \ \ \ \ \ \ \ \ \ $b=6$, $c=8$.}}{ {\rm \ \ \ \ \ \ \ \ \ \ \ \ \ \ \ \ \ \ \ \ \ \ \ $a=5$, $b=7$, $c=6$.}}
\endinsert

\topinsert
\twoline{\mypic{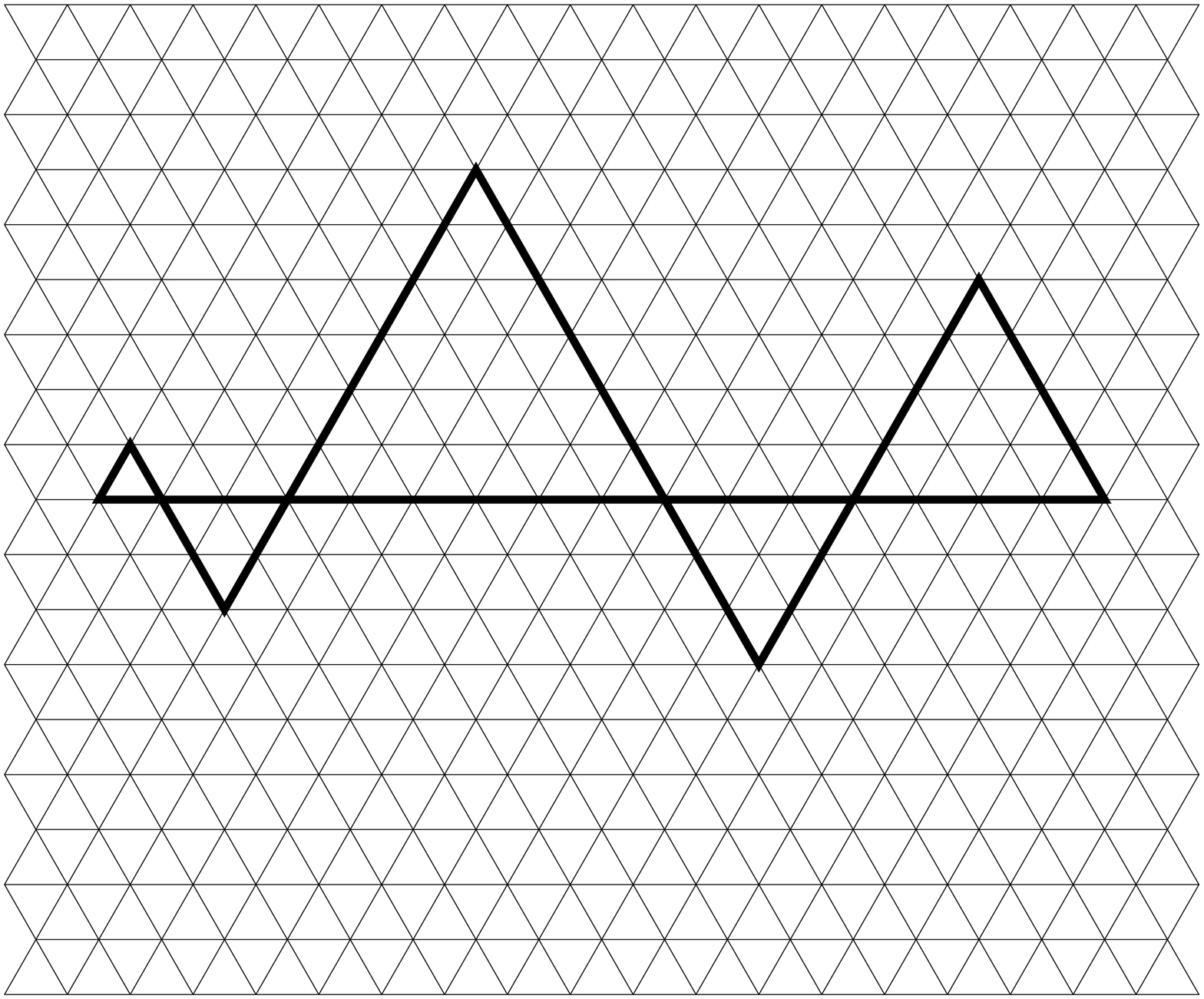}}{\mypic{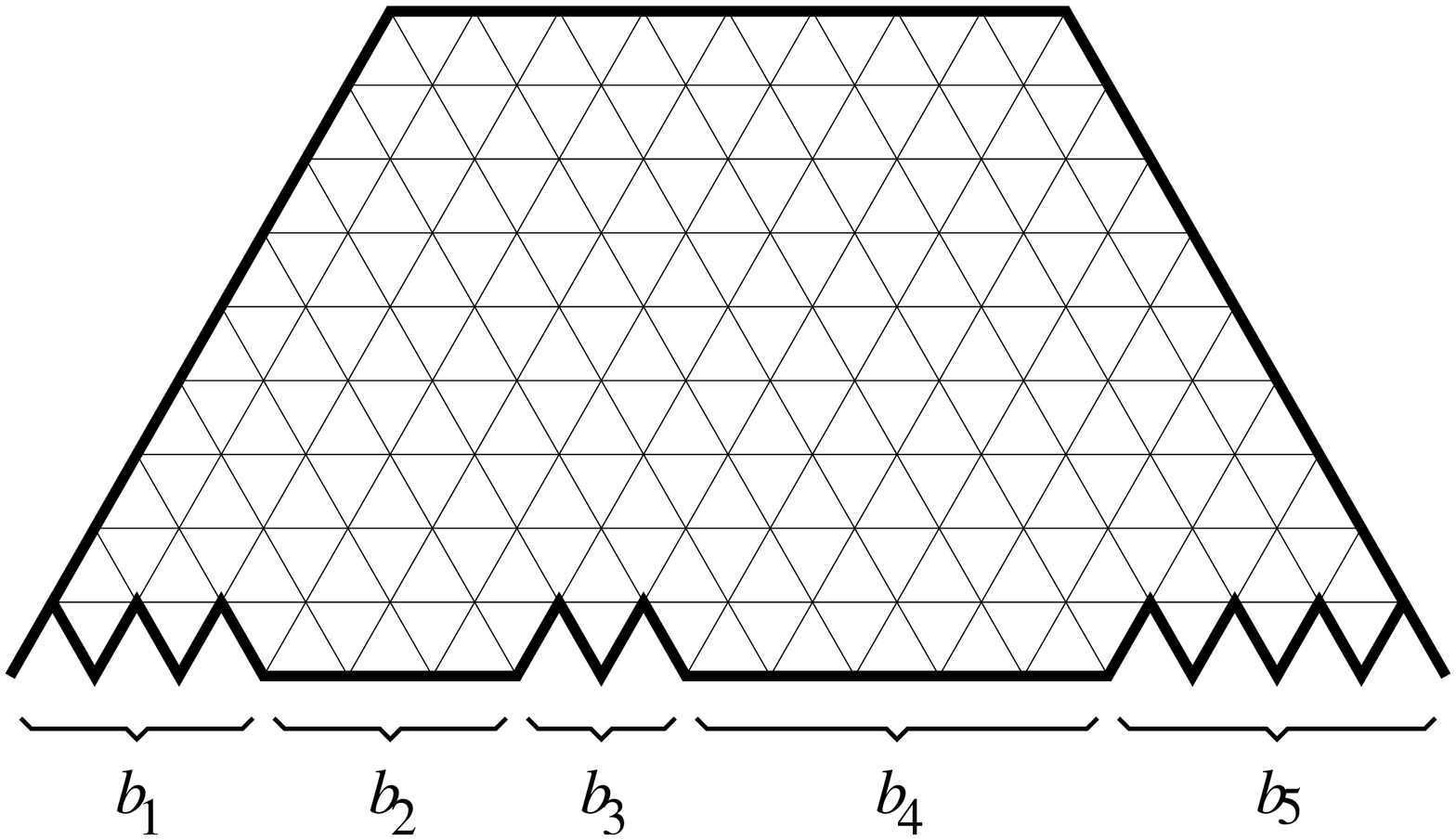}}
\twoline{Figure~{\fac}. {\rm Fern with $a_1=1$, $a_2=2$,}}{\ \ \ \ \ \ Figure~{\fad}. {\rm $S(b_1,b_2,b_3,b_4,b_5)$ for $b_1=3$,}}
\twoline{ {\rm \ \ \ \ \ \ \ \ \ \ \ \ $a_3=6$, $a_4=3$, $a_5=4$.}}{ {\rm \ \ \ \ \ \ \ \ \ \ \ \ \ \ \ \ \ \ \ \ \ \ \ \ \ \ \ \ \ \ $b_2=3$, $b_3=2$, $b_4=5$, $b_5=4$.}}
\endinsert

In this paper we introduce a new structure, called a {\it fern} --- an arbitrary string of triangles of alternating orientations that touch at corners and are lined up along a common axis  (see Figure {\fac} for an example) --- so that the normalized number of tilings of their exterior is given by a product formula in the style of (\eaa). In fact, when the fern has three lobes, one obtaines {\it precisely} formula (\eaa)! Therefore, in some sense, this naturally places MacMahon's formula in a sequence of more general exact enumeration formulas. From this point of view, the situation is more satisfying than the one in \cite{\vf}, as a shamrock only has one lobe in addition to a 3-lobe fern.

Our results concern the exterior of a fern (as its interior has no lozenge tilings, just as it was the case for shamrocks). Let $F(a_1,\dotsc,a_k)$ be the fern whose successive lobes, from left to right, are equilateral triangles of side-lengths $a_1,\dotsc,a_k$, with $a_1$ pointing upwards\footnote{ This is no restriction of generality, as ferns with leftmost lobe pointing downward are obtained by setting $a_1=0$.} (Figure {\fac} illustrates $F(1,2,6,3,4)$). Denote its exterior by $F^*(a_1,\dotsc,a_k)$.

We define the ratio of the number of tilings of the exteriors of the ferns $F(a_1,\dotsc,a_k)$ and $F(a_1+a_3+a_5+\cdots,a_2+a_4+a_6+\cdots)$ as follows. Let $H_N(a_1,\dotsc,a_k)$ be the hexagonal region of side-lengths
alternating between $N+a_1+a_3+a_5+\cdots$ and $N+a_2+a_4+a_6+\cdots$ (the top side being $N+a_2+a_4+a_6+\cdots$), and having the fern $F(a_1,\dotsc,a_k)$ removed from its center (to be precise, $H_N(a_1,\dotsc,a_k)$ is the region $FC_{N,N,N}(a_1,\dotsc,a_k)$ described in the next section). Then we define
$$
\spreadlines{3\jot}
\align
&
\frac{\M(F^*(a,\dotsc,a_k))}{\M(F^*(a_1+a_3+a_5+\cdots,a_2+a_4+a_6+\cdots))}
:=
\\
&\ \ \ \ \ \ \ \ \ \ \ \ 
\lim_{N\to\infty}\frac{\M(H_N(a_1,\dotsc,a_k))}{\M(H_N(a_1+a_3+a_5+\cdots,a_2+a_4+a_6+\cdots))},
\tag\eac
\endalign
$$
where for a lattice region $R$, $\M(R)$ denotes the number of lozenge tilings of $R$.

%For integers $0<m\leq n$ and $1\leq x_1<\cdots< x_m\leq n$, denote by $H(x_1,\dotsc,x_m)$ the semihexagon of base $n$

%The following formula is Cohn, Larsen and Propp's \cite{\CLP} translation to lozenge tilings of a classical result of Gelfand and Tsetlin \cite{\GT}.

%\proclaim{Proposition \tcb} Let $T_{m,n}(x_1,\dotsc,x_n)$ be the region obtained from the trapezoid of side lengths $m$, $n$, $m+n$, $n$ (clockwise from bottom) by removing the down-pointing unit triangles from along its top that are in positions $x_1,x_2,\dotsc,x_n$ as counted from left to right\footnote{ For short, we sometimes refer to a region of this kind as a $T$-region, or a $T$-type region.}. Then
%$$
%\M(T_{m,n}(x_1,\dotsc,x_n))=\prod_{1\leq i<j\leq n}\frac{x_j-x_i}{j-i}.
%\tag\ecf
%$$

%\endproclaim  

Let $s(b_1,b_2,\dotsc,b_l)$ denote the number of lozenge tilings of the semihexagonal region $S(b_1,b_2,\dotsc,b_l)$ with the leftmost $b_1$ up-pointing unit triangles on its base removed, the next segment of length $b_2$ intact, the following $b_3$ removed, and so on (see an illustration in Figure {\fad}; note that the $b_i$'s determine the lengths of all four sides of the semihexagon). By the Cohn-Larsen-Propp \cite{\CLP} interpretation of the Gelfand-Tsetlin result \cite{\GT} we have that\footnote{ The first equality in (\ead) holds due to forced lozenges in the tilings of $S(b_1,b_2,\dotsc,b_{2l})$, after whose removal one is left precisely with the region $S(b_1,b_2,\dotsc,b_{2l-1})$.}${}^{,}$\footnote{ We include here the original formula for convenience. Let $T_{m,n}(x_1,\dotsc,x_n)$ be the region obtained from the trapezoid of side lengths $m$, $n$, $m+n$, $n$ (clockwise from top) by removing the up-pointing unit triangles from along its bottom that are in positions $x_1,x_2,\dotsc,x_n$ as counted from left to right. Then
$$
\M(T_{m,n}(x_1,\dotsc,x_n))=\prod_{1\leq i<j\leq n}\frac{x_j-x_i}{j-i}.
$$}
$$
s(b_1,b_2,\dotsc,b_{2l})=s(b_1,b_2,\dotsc,b_{2l-1})=
\frac{\prod_{\text{\rm $1\leq i<j\leq 2l-1$, $j-i+1$ odd}}\h(b_i+b_{i+1}+\cdots+b_j)}{\prod_{\text{\rm $1\leq i<j\leq 2l-1$, $j-i+1$ even}}\h(b_i+b_{i+1}+\cdots+b_j)}.
\tag\ead
$$

The dual MacMahon theorem we obtain in this paper is the following.

\proclaim{Theorem \taa} For any non-negative integers $a_1,\dotsc,a_k$ we have
$$
\spreadlines{3\jot}
\align
&
\frac{\M(F^*(a_1,\dotsc,a_k))}{\M(F^*(a_1+a_3+a_5+\cdots,a_2+a_4+a_6+\cdots))}
=
\\
&\ \ \ \ \ \ \ \ \ \ \ \ \ \ \ \ \ \ \ \ \ \ \ \ \ \ \ \ \ \ \ \ \ \ \ \ \ \ \ 
s(a_1,a_2,\dotsc,a_{k-1})\,s(a_2,a_3,\dotsc,a_k).
\tag\eae
\endalign
$$

\endproclaim

In the special case when $k=3$, this becomes
$$                                                                                              
%\spreadlines{3\jot}                                                                             
%\align                                                                                          
\frac{\M(F^*(a_1,a_2,a_3))}{\M(F^*(a_1+a_3,a_2))}       
%&
=                                                                                              
\frac{\h(a_1)\h(a_2)\h(a_3)\h(a_1+a_2+a_3)}{\h(a_1+a_2)\h(a_1+a_3)\h(a_2+a_3)}
%\\                                                                                              
%&
=P(a_1,a_2,a_3),
\tag\eaf
%\endalign
$$
giving thus prescisely MacMahon's expression.

Formula (\eae) is illustrated geometrically in Figure {\fea} (see section 5). Theorem {\taa} will 
follow as a consequence of a more general result (see Theorem {\tba}), which we describe in the next section.

We end this section by a brief account on how the exact formulas in \cite{\vf} and the current paper were found.
%, and how the proofs were arrived at. 
The author noticed several years ago, by working out concrete examples, that the correlation\footnote{We denote by $\omega(R)$ the correlation of the region $R$ on the triangular lattice as it is defined in our earlier work \cite{\ec} and \cite{\ov}.} $\omega(S)$ of shamrocks, as well as the correlation $\omega(F)$ of ferns, seemed to be gived by simple product formulas. He later noticed that there is a way to enclose each of these two structures near the center of a hexagon so that the number of tilings of the resulting hexagon with the corresponding hole appears to be given by simple product formulas. 

The case of the shamrock was proved in \cite{\vf}. 
%The initial approach was to use the identification of factors method, which proved successful in a series of previous articles, but in the face of mounting technical difficulties, the author noticed that the proof follows using Kuo's graphical condensation method \cite{\Kuo}. 
The main purpose of this paper is to prove the remaining case of the fern\footnote{It may amuse the reader to know that initially the exact formulas in \cite{\vf} and the current paper were noticed in the case when there are just two lobes, when both structures become shaped like a bowtie, or farfalle. When the author extended them to the general form, he was visiting the University of Vienna, and found it natural to call the four-lobed and many-lobed structures vierfalle and farviele, respectively.}.

The family consisting of shamrocks and  ferns has the following simple description: It consists of connected unions of triangles touching only at vertices in which, if regarded as an island in a surrounding sea, 
%(in a sea of dimers?), 
all gulfs are open, in the sense that no gulf has two parallel sides (i.e., no ``fjords''; see this in Figures {\fab} and {\fac}). Data suggests that these are the only structures that lead to simple product formulas of this kind. Any insight into why this happens would be very interesting.

\mysec{2. Precise statement of results}

% Better: Describe regions, then point put that for k=1 we get cored hex

For non-negative integers $a_1,\dotsc,a_k$, define the {\it fern} $F(a_1,\dotsc,a_k)$ to be a string of $k$ lattice triangles lined up along along a horizontal lattice line, touching at their vertices, alternately oriented up and down and having sizes $a_1,\dotsc,a_k$ as encountered from left to right (with the leftmost oriented up; see footnote 1); the black structure in Figure {\fba} is a 4-lobed fern with lobes of sizes 1, 2, 6 and 2. The {\it base} (or basepoint) of a fern is its leftmost point (the left vertex of the first lobe); for the fern in Figure {\fba}, the base is indicated by a black dot.

%The regions we will be concerned with in this paper are defined as follows. Let $H$ be the hexagon of side-lengths $x+a_2+a_4+a_6+\cdots$, $y+a_1+a_3+a_5+\cdots$, $z+a_2+a_4+a_6+\cdots$, $x+a_1+a_3+a_5+\cdots$, $y+a_2+a_4+a_6+\cdots$, $z+a_1+a_3+a_5+\cdots$ (clockwise from top). Then our regions are obtained from $H$ by taking out from near the center of $H$ the fern $F(a_1,\dotsc,a_k)$. The details depend on the relative parities of $x$, $y$ and $z$, and lead us to four versions of our regions. We call the hexagon $\bar{H}$ of side-lengths $x$, $y$, $z$, $x$, $y$, $z$ (clockwise from top) that fits in the western corner of $H$ the {\it auxilliary hexagon} (in Figures {\fba} and~{\fbb} the auxilliary hexagon is indicated by a dotted line). The details are as follows.

The regions we will be concerned with in this paper are defined as follows. Set
$$
\align
o&:=a_1+a_3+a_5+\cdots\tag\ebaa\\
e&:=a_2+a_4+a_6+\cdots\tag\ebab
\endalign
$$
and let $H$ be the hexagon of side-lengths $x+e$, $y+o$, $z+e$, $x+o$, $y+e$, $z+o$ (clockwise from top). We call the hexagon $\bar{H}$ of side-lengths $x$, $y$, $z$, $x$, $y$, $z$ (clockwise from top) that fits in the western corner of $H$ the {\it auxilliary hexagon} (in Figures {\fba} and~{\fbb} the auxilliary hexagon is indicated by a thick dotted line). Our regions are obtained by taking out from $H$ a certain translation of the fern that places its basepoint either at the very center of the auxilliary hexagon (when this center is a lattice point), or as close to it as possible. The details, which depend on the relative parities of $x$, $y$ and $z$, are as follows.

\topinsert
\twoline{\mypic{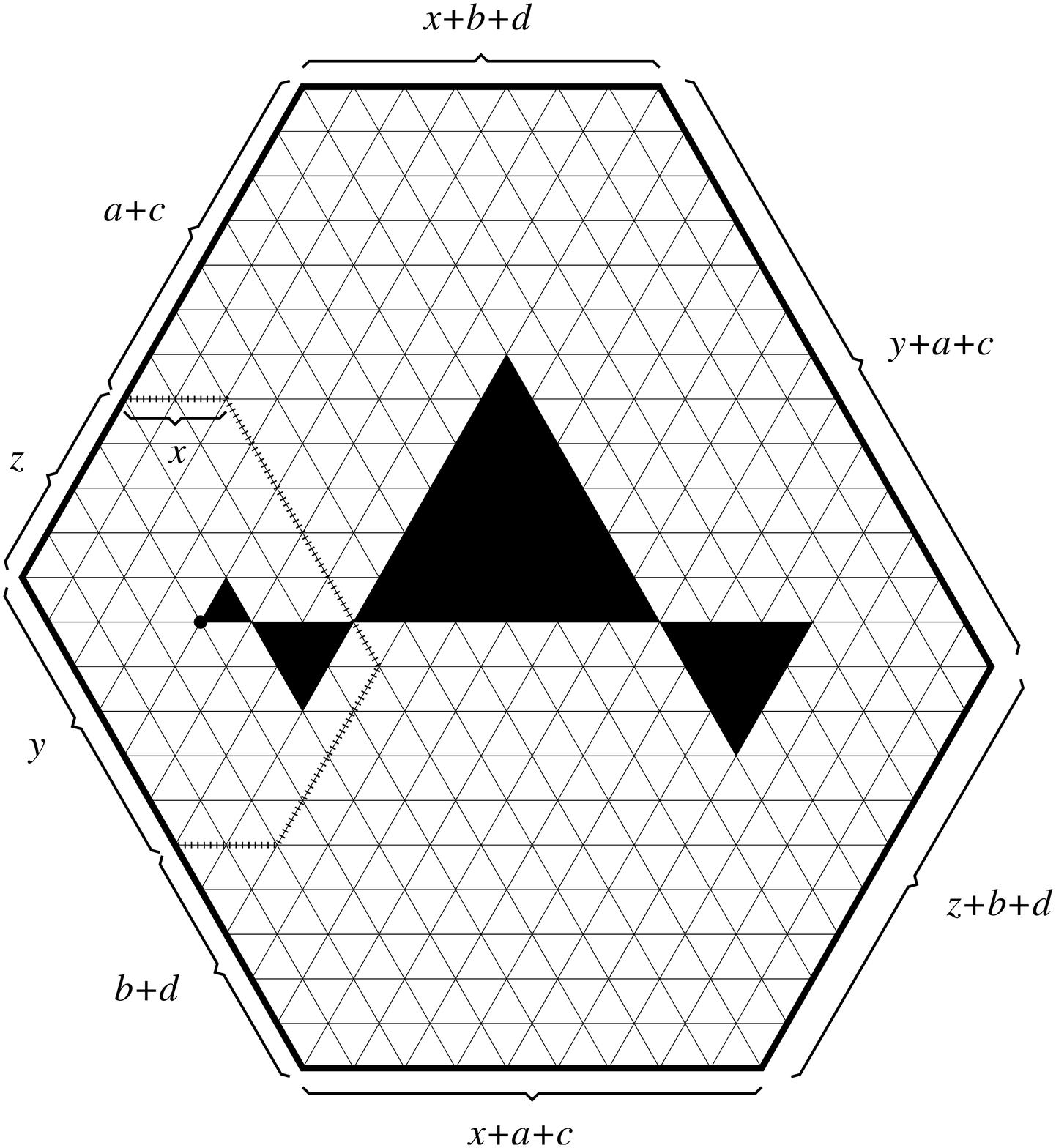}}{\mypic{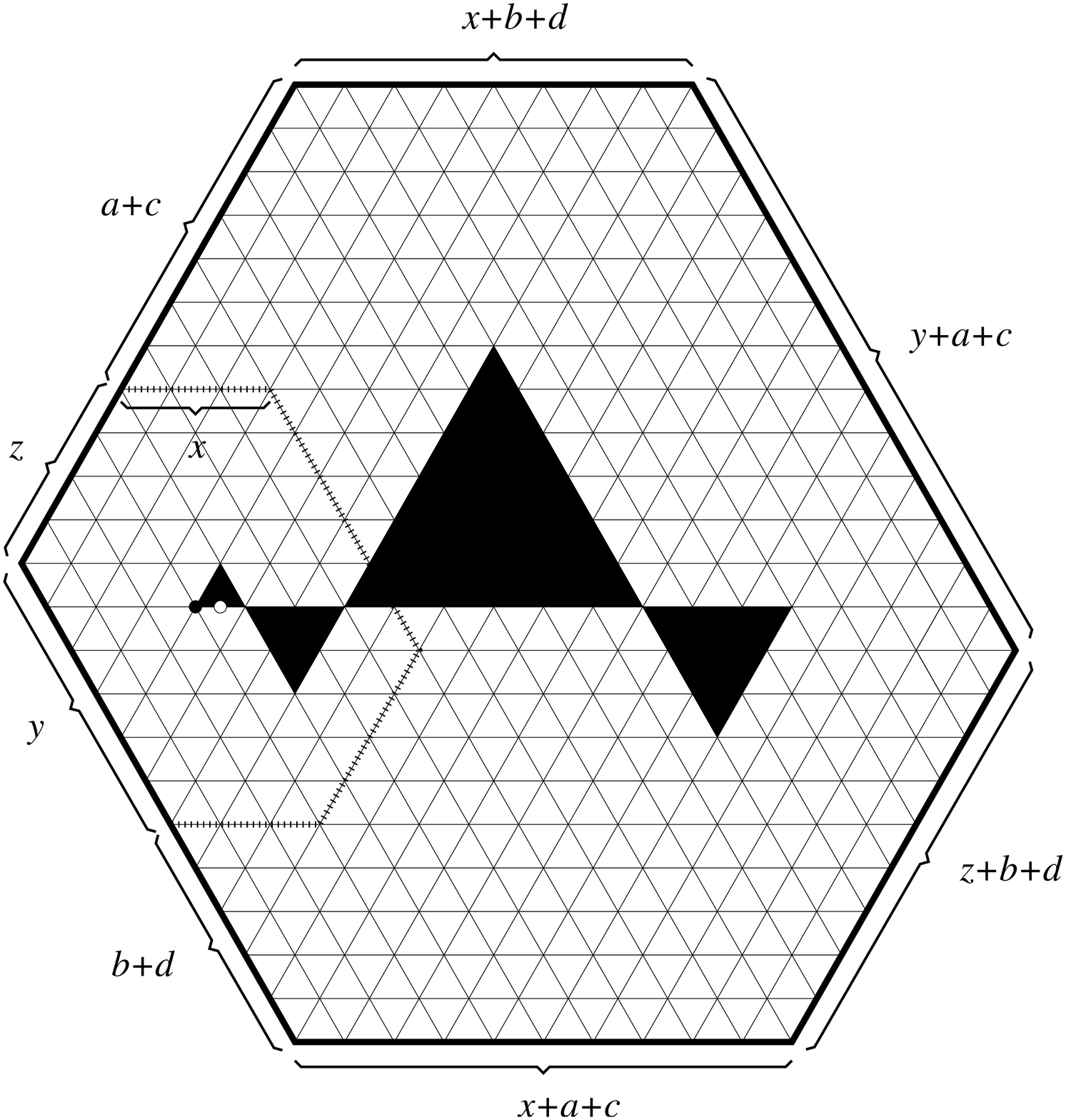}}
\medskip
\centerline{Figure~2.1. {\rm The $F$-cored hexagons $FC^{\bigodot}_{2,6,4}(1,2,6,3)$ (left) and $FC^{\leftarrow}_{3,6,4}(1,2,6,3)$ (right).}}
%\endinsert
\medskip
%\topinsert
\twoline{\mypic{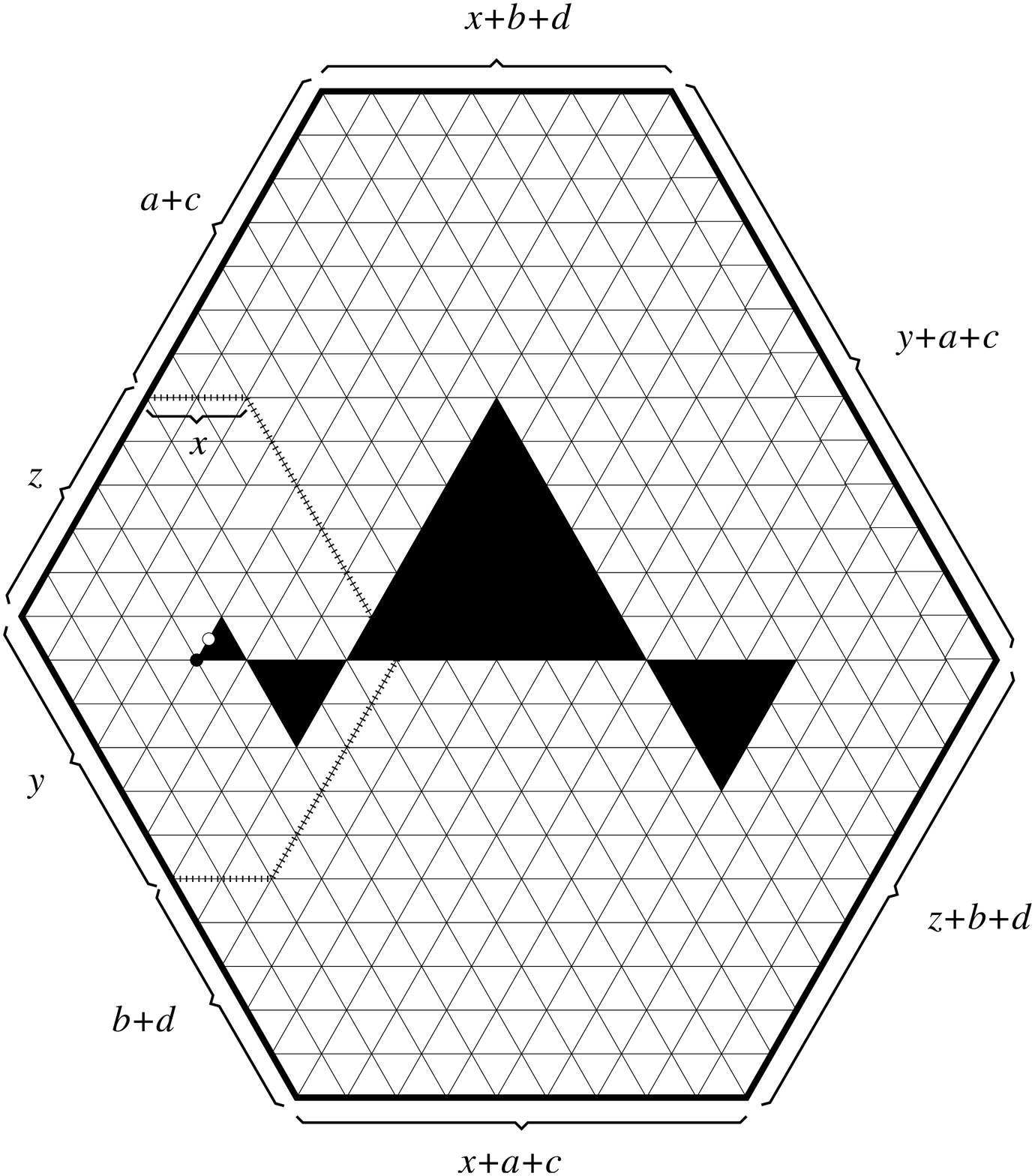}}{\mypic{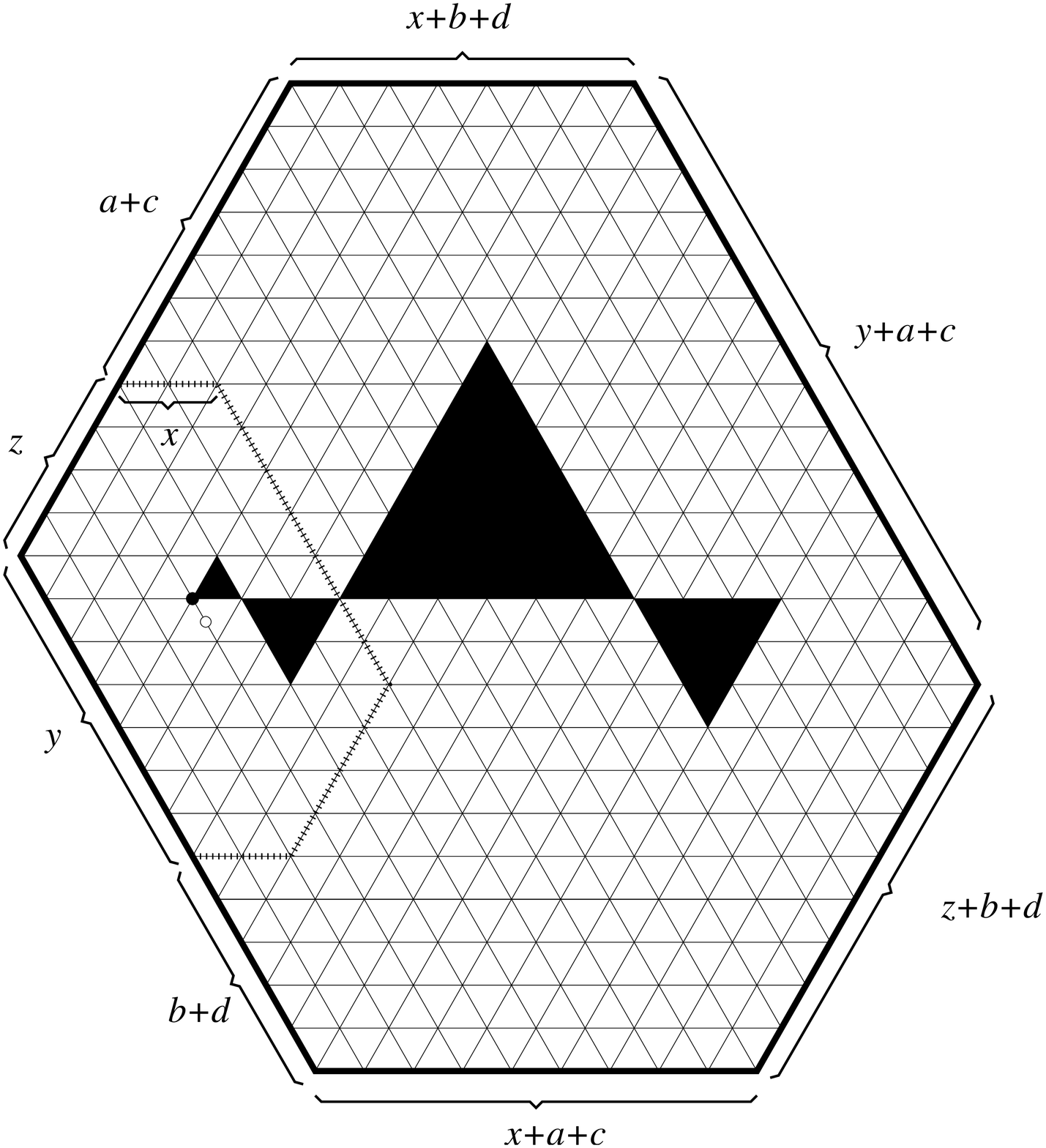}}
\medskip
\centerline{Figure~2.2. {\rm The $F$-cored hexagons $FC^{\swarrow}_{2,6,5}(1,2,6,3)$ (left) and $FC^{\nwarrow}_{2,7,4}(1,2,6,3)$ (right).}}
\endinsert

The center of $\bar{H}$ is a lattice point precisely when $x$, $y$ and $z$ have the same parity. If this is the case, place the fern $F(a_1,\dotsc,a_k)$ inside the hexagon $H$ so that the base of the fern 
%(i.e., its leftmost point) is $(x+y)/2$ lattice spacings from the northwestern side of $H$ and $(x+z)/2$ lattice spacings from the southwestern side of $H$ --- i.e., 
is at the center of the auxilliary hexagon $\bar{H}$ (an example is illustrated in Figure~{\fba}, left).
%; the base is marked by a black dot). 
Denote by $FC^{\bigodot}_{x,y,z}(a_1,\dotsc,a_k)$ the region obtained from $H$ by removing the fern placed in this position (the superscript is a visual reminder of the fact that the fern's base is precisely at the center of the auxilliary hexagon).

If $x$ has opposite parity to $y$ and $z$, draw the fern $F(a_1,\dotsc,a_k)$ inside $H$ so that the base of the fern 
%is $(x+y-1)/2$ lattice spacings from the northwestern side of $H$ and $(x+z-1)/2$ lattice spacings from the southwestern side of $H$. 
is half a lattice spacing to the west of the center of the auxilliary hexagon (this center is marked by a white dot in Figure {\fba}, right).
Denote by $FC^\leftarrow_{x,y,z}(a_1,\dotsc,a_k)$ the region obtained by taking out from $H$ this placement of the fern; the arrow at the superscript records the relative position of the base versus the center of the auxilliary hexagon.
%records the fact that the base of the fern is half a lattice spacing to the west of the center of the auxilliary hexagon; the latter is marked by a white dot in the figure).

If $z$ has parity opposite to $x$ and $y$, we place the fern so that its base is half a lattice spacing to the southwest of the center of $\bar{H}$
and denote the resulting region by $FC^\swarrow_{x,y,z}(a_1,\dotsc,a_k)$ (see Figure {\fbb}, left for an illustration). Finally, if $y$ has parity opposite to $x$ and $z$, place the fern so that its base is half a lattice spacing to the the northwest of the center of $\bar{H}$, and denote the resulting region by $FC^\nwarrow_{x,y,z}(a_1,\dotsc,a_k)$ (see Figure {\fbb}, right for an illustration).

It will be convenient to have a common notation for these four families of regions. We define the {\it $F$-cored hexagon} $FC_{x,y,z}(a_1,\dotsc,a_k)$ by
$$
FC_{x,y,z}(a_1,\dotsc,a_k):=
\cases
FC^{\bigodot}_{x,y,z}(a_1,\dotsc,a_k), & {\text{\rm if $x$, $y$, $z$ same parity}} \\
FC^\leftarrow_{x,y,z}(a_1,\dotsc,a_k), & {\text{\rm if $x$ has parity opposite to $y$, $z$}} \\
FC^\swarrow_{x,y,z}(a_1,\dotsc,a_k), & {\text{\rm if $z$ has parity opposite to $x$, $y$}} \\
FC^\nwarrow_{x,y,z}(a_1,\dotsc,a_k), & {\text{\rm if $y$ has parity opposite to $x$, $z$}}
\endcases
\tag\eba
$$
For $k=1$, these become the cored hexagons we studied in \cite{\cekz}. 

\smallpagebreak
\flushpar
{\bf Remark 1.} One may wonder why is it that we do not consider three more types of such regions, with the base of the fern displaced half a lattice spacing from the center of the auxilliary hexagon in the eastern, northeastern, or southeastern directions. The reason for this is that the latter can be viewed as $180^\circ$-degree rotations of the last three families in (\eba) (this is visually apparent when $k$ is even; for odd $k$, after the rotation the leftmost lobe points down instead of up, but then we regard the fern as having $k+1$ lobes, the first of which is of side length 0).

By contrast, in \cite{\vf} (as well as in \cite{\cekz}) it was enough to consider one of the six possible off-center placements of the shamrock in the hexagon, because all six can be obtained from any of them by rotations by $120^\circ$ and reflection across the vertical, symmetries that preserve both the structure of the hexagon and of the shamrock (but not of the fern).

\medskip
The main result of this paper is the following.

\proclaim{Theorem \tba} Let $x$, $y$, $z$ and $a_1,\dotsc,a_k$ be non-negative integers. Then the number of lozenge tilings of the $F$-cored hexagon $FC_{a,b,c}(a_1,\dotsc,a_k)$ is given by
%\footnote{ For brevity, we denote by $a_1+a_3+\cdots$ the sum of the odd-indexed $a_i$'s, i.e. the sum of the side lengths of the up-pointing lobes in the fern; $a_2+a_4+\cdots$ has the analogous meaning.}
$$
\spreadlines{3\jot}
\align
&
\frac
{\M(FC_{x,y,z}(a_1,\dotsc,a_k))}
{\M(FC_{x,y,z}(a_1+a_3+a_5+\cdots,a_2+a_4+a_6+\cdots))}
=
s(a_1,\dotsc,a_{k-1})s(a_2,\dotsc,a_k)
\\
&\ \ \ \ \ \ \ \ \ \ \ \ \ \ \ \ \ \ \ \ \ \ \ 
\times
\frac
{\h(\lfloor \frac{x+z}{2}\rfloor+a_1+a_3+\cdots)}
{\h(\lfloor \frac{x+y}{2}\rfloor+a_1+a_3+\cdots)}
\frac
{\h(\lceil \frac{x+z}{2}\rceil+a_2+a_4+\cdots)}
{\h(\lceil \frac{x+y}{2}\rceil+a_2+a_4+\cdots)}
\\
&\ \ \ \ \ \ \ \ \ \ 
\times
\prod_{1\leq 2i+1\leq k} 
\frac
{\h(\lfloor \frac{x+y}{2}\rfloor+a_1+\cdots+a_{2i+1})}
{\h(\lfloor \frac{x+z}{2}\rfloor+a_1+\cdots+a_{2i+1})}
\frac
{\h(\lceil \frac{x+y}{2}\rceil+\overline{a_1+\cdots+a_{2i+1}})}
{\h(\lceil \frac{x+z}{2}\rceil+\overline{a_1+\cdots+a_{2i+1}})}
\\
&\ \ \ \ \ \ \ \ \ \ \ 
\times
\prod_{1< 2i< k} 
\frac
{\h(\lfloor \frac{x+z}{2}\rfloor+a_1+\cdots+a_{2i})}
{\h(\lfloor \frac{x+y}{2}\rfloor+a_1+\cdots+a_{2i})}
\frac
{\h(\lceil \frac{x+z}{2}\rceil+\overline{a_1+\cdots+a_{2i}})}
{\h(\lceil \frac{x+y}{2}\rceil+\overline{a_1+\cdots+a_{2i}})}
,
\tag\ebb
\endalign
$$
where $H$ and $s$ are defined by (1.2) and (1.4), and $\overline{a_1+\cdots+a_i}$ stands for $a_{i+1}+\cdots+a_k$.
\endproclaim

We note that when the fern has only two lobes, $F$-cored hexagons are specializations of the $S$-cored hexagons $SC_{x,y,z}(a,b,c,m)$ whose lozenge tilings were enumerated in \cite{\vf}. More precisely, we have that 
$$
\spreadlines{3\jot}
\align
FC^{\bigodot}_{x,y,z}(a,b)&=SC_{x,y,z}(0,0,b,a)\tag\ebc
\\
FC^\leftarrow_{x,y,z}(a,b)&=SC_{x,y,z}(0,0,b,a)\tag\ebd
\\
FC^\swarrow_{x,y,z}(a,b)&=SC_{z,y,x}(a,0,0,b)\tag\ebe
\\
FC^\nwarrow_{x,y,z}(a,b)&=SC_{y,z,x}(0,a,0,b).\tag\ebf
\endalign
$$
With the product formulas for the number of lozenge tilings of $S$-cored hexagons given in \cite{\vf, Theorems\,2.1,\,2.2}, we see that (\ebb) does indeed provide explicit expressions for the number of lozenge tilings of $F$-cored hexagons with an arbitrary number of lobes. This is made explicit at the end of Remark 2 below.

\medskip
\flushpar
{\bf Remark 2.} It turns out that the somewhat lengthy product formulas supplied by \cite{\vf, Theorems\,2.1,\,2.2} can be simplified significantly for the case when the shamrock has only one lobe (i.e., the shamrock becomes a fern with two lobes, which is the case needed in (\ebc)--(\ebf)). Indeed, we get that if $y$ and $z$ have the same parity\footnote{ 
%This case covers the general situation, as the case when $x$ and $y$ or $x$ and $z$ have the same parity can be reduced to it by a permutation of $x$, $y$ and $z$, and by viewing the removed 2-lobe fern as a shamrock in a different way. 
We are using here the fact --- which went unnoticed at the time of writing of \cite{\vf} --- that the expression in \cite{\vf, Theorem 2.2} specializes to the one in \cite{\vf, Theorem 2.1} when $x$, $y$ and $z$ have the same parity. Thus Theorems 2.1 and 2.2 of \cite{\vf} could have been stated as a single result, holding when $y$ and $z$ have the same parity.}
$$
\spreadlines{3\jot}
\align
&
\frac
{\M(FC_{x,y,z}(a,b))}
{\M(C_{x,y,z}(a+b))}
=
\frac{\h(a)\h(b)}{\h(a+b)}
\frac{\h(\lceil\frac{x+z}{2}\rceil)\h(\frac{y+z}{2})}
{\h(\lceil\frac{x+y}{2}\rceil)}
%\\
%&
%\times
\frac{\h(a+\lfloor\frac{x+y}{2}\rfloor)\h(b+\lceil\frac{x+y}{2}\rceil)}
{\h(a+b+\lfloor\frac{x+y}{2}\rfloor)}
\\
&\ \ \ \ \ \ \ \ \ \ \ \ \ \ \ \ \ \ \ \ 
\times
\frac{\h(a+b+\lfloor\frac{x+z}{2}\rfloor)}
{\h(a+\lfloor\frac{x+z}{2}\rfloor)\h(b+\lceil\frac{x+z}{2}\rceil)}
\frac{\h(a+b+\frac{y+z}{2})}
{\h(a+\frac{y+z}{2})\h(b+\frac{y+z}{2})}
\tag\ebg
\endalign
$$
(the regions $C_{x,y,z}(m)$ are the cored hexagons considered in \cite{\cekz}, which are just the special case $SC_{x,y,z}(0,0,0,m)$ of the $S$-cored hexagons of \cite{\vf}). According to Theorems 1 and 2 of \cite{\cekz}, if $y$ and $z$ have the same parity, the number of lozenge tilings of the cored hexagon is given by\footnote{ Here we are using that, as in the previous footnote, Theorems 1 and 2 of \cite{\cekz} can be combined into the single case of $y$ and $z$ having the same parity.}
$$
\spreadlines{3\jot}
\align
&
\M(C_{x,y,z}(m))=
\\
&\ \ \ \ \ \ 
\frac
{\h(x+m)\h(y+m)\h(z+m)\h(x+y+z+m)
\h(\lfloor\frac{x+y+z}{2}\rfloor+m)\h(\lceil\frac{x+y+z}{2}\rceil+m)}
{\h(x+y+m)\h(x+z+m)\h(y+z+m)
\h(\lceil\frac{x+y}{2}\rceil+m)\h(\lfloor\frac{x+z}{2}\rfloor+m)
\h(\frac{y+z}{2}+m)}
\\
&\ \ \,
\times
\frac
{\h(\frac{m}{2})^2
\h(\lfloor\frac{x}{2}\rfloor)\h(\lceil\frac{x}{2}\rceil)
\h(\lfloor\frac{y}{2}\rfloor)\h(\lceil\frac{y}{2}\rceil)
\h(\lfloor\frac{z}{2}\rfloor)\h(\lceil\frac{z}{2}\rceil)}
{\h(\lfloor\frac{x}{2}\rfloor+\frac{m}{2})\h(\lceil\frac{x}{2}\rceil+\frac{m}{2})
\h(\lfloor\frac{y}{2}\rfloor+\frac{m}{2})\h(\lceil\frac{y}{2}\rceil+\frac{m}{2})
\h(\lfloor\frac{z}{2}\rfloor+\frac{m}{2})\h(\lceil\frac{z}{2}\rceil+\frac{m}{2})}
\\
&\ \ \,
\times
\frac
{\h(\lfloor\frac{x+y}{2}\rfloor+\frac{m}{2})
\h(\lceil\frac{x+y}{2}\rceil+\frac{m}{2})
\h(\lfloor\frac{x+z}{2}\rfloor+\frac{m}{2})
\h(\lceil\frac{x+z}{2}\rceil+\frac{m}{2})
\h(\frac{y+z}{2})^2}
{\h(\lfloor\frac{x+y+z}{2}\rfloor+\frac{m}{2})
\h(\lceil\frac{x+y+z}{2}\rceil+\frac{m}{2})
\h(\lfloor\frac{x+y}{2}\rfloor)\h(\lceil\frac{x+z}{2}\rceil)\h(\frac{y+z}{2})}
.
\tag\ebh
\endalign
$$
Then the explicit expression for the number of tilings of $F$-cored hexagons is obtained by formulas (\ebb), (\ebg) and (\ebh).

The remaining parity cases are given by analogous formulas. For completeness we include them below.

When $x$ and $y$ have the same parity, the analogs of (\ebg) and (\ebh) are
$$
\spreadlines{3\jot}
\align
&
\frac
{\M(FC_{x,y,z}(a,b))}
{\M(C_{x,y,z}(a+b))}
=
\frac{\h(a)\h(b)}{\h(a+b)}
\frac{\h(\lfloor\frac{x+z}{2}\rfloor)\h(\lceil\frac{y+z}{2}\rceil)}
{\h(\frac{x+y}{2})}
%\\
%&
%\times
\frac{\h(a+\frac{x+y}{2})\h(b+\frac{x+y}{2})}
{\h(a+b+\frac{x+y}{2})}
\\
&\ \ \ \ \ \ \ \ \ \ \ \ \ \ \ \ \ \ \ \ 
\times
\frac{\h(a+b+\lceil\frac{x+z}{2}\rceil)}
{\h(a+\lceil\frac{x+z}{2}\rceil)\h(b+\lfloor\frac{x+z}{2}\rfloor)}
\frac{\h(a+b+\lfloor\frac{y+z}{2}\rfloor)}
{\h(a+\lfloor\frac{y+z}{2}\rfloor)\h(b+\lceil\frac{y+z}{2}\rceil)}
\tag\ebga
\endalign
$$
and
$$
\M(C_{x,y,z}(m))=
\text{\rm right hand side of (\ebh)}|_{x\leftarrow z,y\leftarrow x,z\leftarrow y}
.
\tag\ebha
$$
Finally, if $x$ and $z$ have the same parity, the analogs of (\ebg) and (\ebh) are
$$
\spreadlines{3\jot}
\align
&
\frac
{\M(FC_{x,y,z}(a,b))}
{\M(C_{x,y,z}(a+b))}
=
\frac{\h(a)\h(b)}{\h(a+b)}
\frac{\h(\frac{x+z}{2})\h(\lfloor\frac{y+z}{2}\rfloor)}
{\h(\lfloor\frac{x+y}{2}\rfloor)}
%\\
%&
%\times
\frac{\h(a+\lceil\frac{x+y}{2}\rceil)\h(b+\lfloor\frac{x+y}{2}\rfloor)}
{\h(a+b+\lceil\frac{x+y}{2}\rceil)}
\\
&\ \ \ \ \ \ \ \ \ \ \ \ \ \ \ \ \ \ \ \ 
\times
\frac{\h(a+b+\frac{x+z}{2})}
{\h(a+\frac{x+z}{2})\h(b+\frac{x+z}{2})}
\frac{\h(a+b+\lceil\frac{y+z}{2}\rceil)}
{\h(a+\lceil\frac{y+z}{2}\rceil)\h(b+\lfloor\frac{y+z}{2}\rfloor)}
\tag\ebgb
\endalign
$$
and
$$
\M(C_{x,y,z}(m))=
\text{\rm right hand side of (\ebh)}|_{x\leftarrow y,y\leftarrow z,z\leftarrow x}
.
\tag\ebhb
$$

\medskip
\flushpar
{\bf Remark 3.} Compared to the apparent lack of structure and the lenghtiness of the expressions in \cite{\vf} for the $S$-cored hexagons --- each of which takes up nearly a full page 
%(although see footnote 7) 
--- the formula for the $F$-cored hexagons given by (\ebb) and (\ebg)--(\ebh) (or (\ebga)--(\ebha), resp. (\ebgb)--(\ebhb)) is remarkably structured and brief. 

%Note that in [cekz], [vf] there are two cases. Actually those are covered by 2nd theorems there! Our formula (2.2) (and (2.7)) is much more simple elegant!

\medskip

Theorem {\tba} generalizes the main results of \cite{\cekz}, by introducing an arbitrary number of new parameters to the geometry of the core (the sizes of the fern's lobes; a different, 3-parameter generalization of the results of \cite{\cekz}  was given in \cite{\vf}). This results in a new, multi-parameter generalization of MacMahon's theorem (\eaa).

%\bigskip
%\bigskip
%\bigskip
%\bigskip
%\bigskip

\mysec{3. Proof of Theorem 2.1}

\nopagebreak

% The base cases x=0, or y=0, or z=0 reduce to S-regions like in the introduction!

\topinsert

\twoline{\mypic{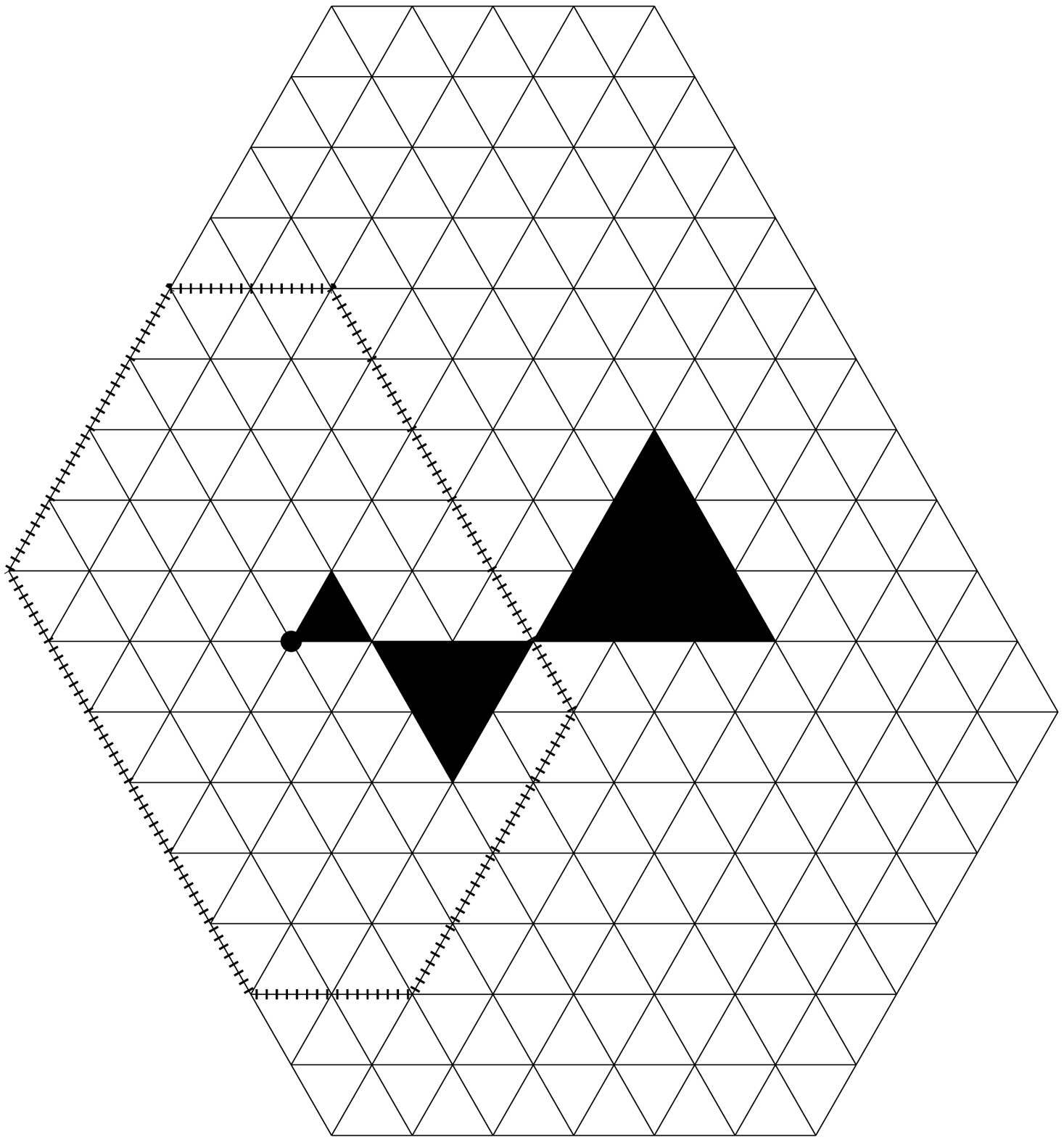}}{\mypic{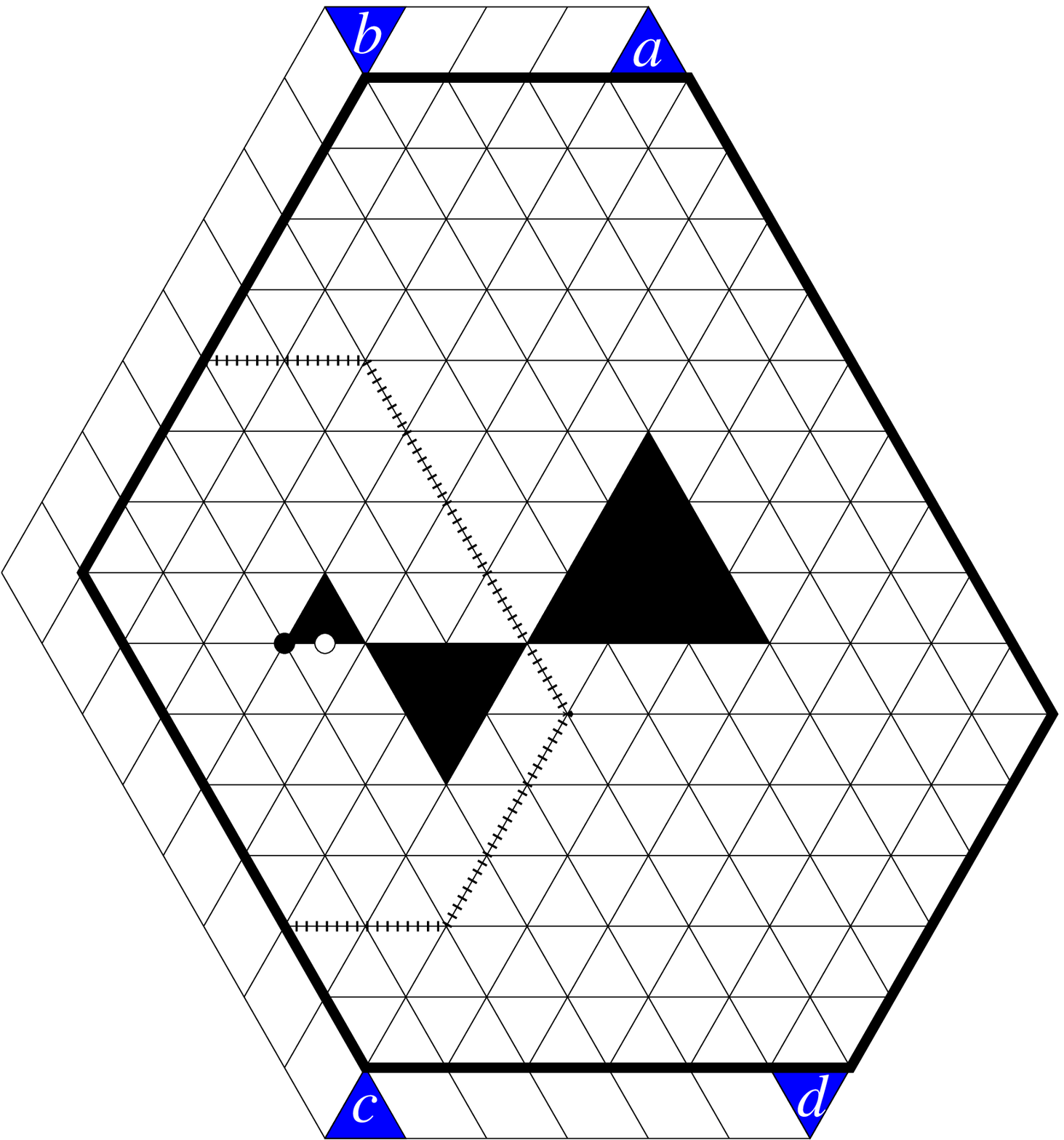}}
\bigskip

\twoline{\mypic{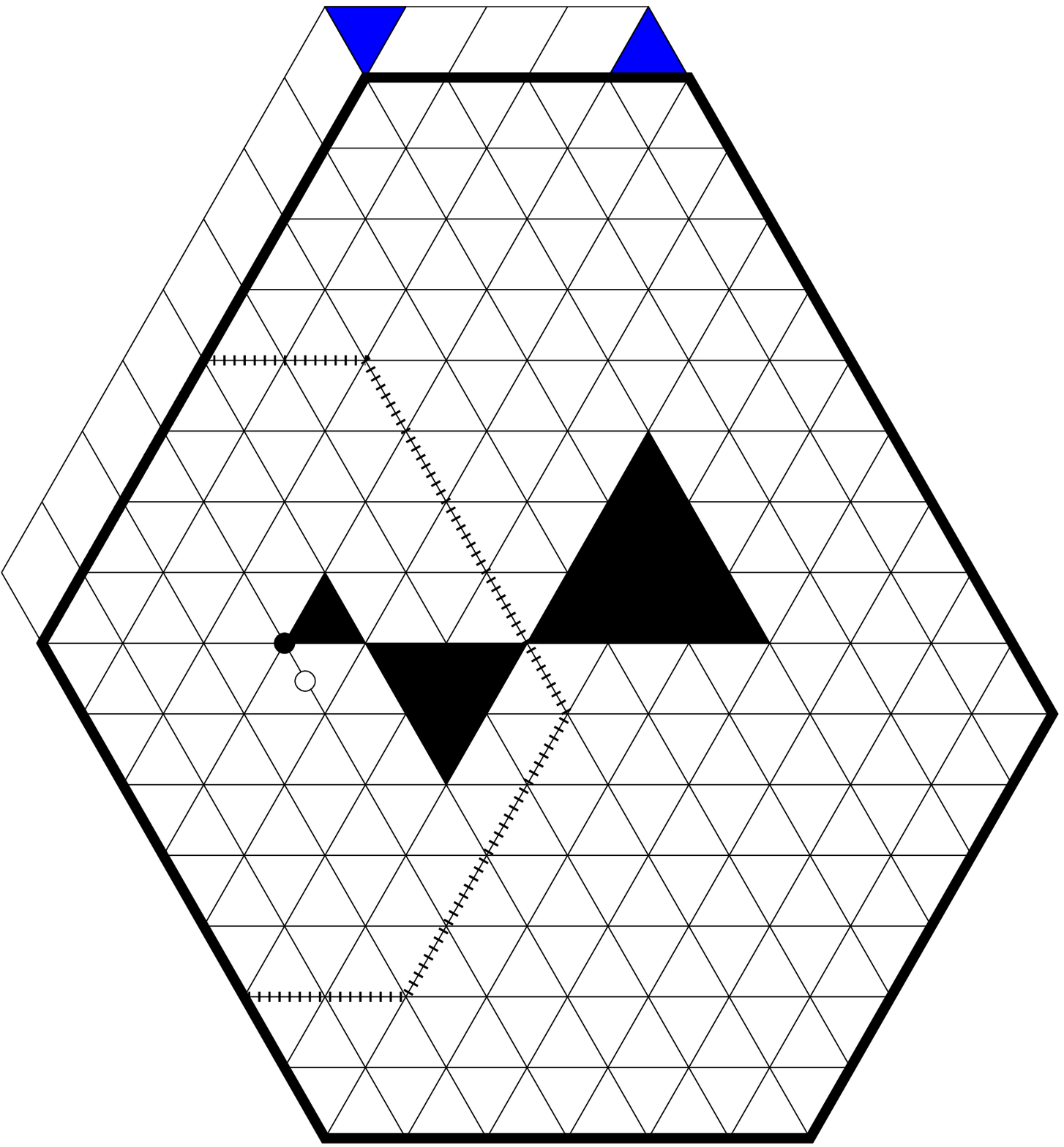}}{\mypic{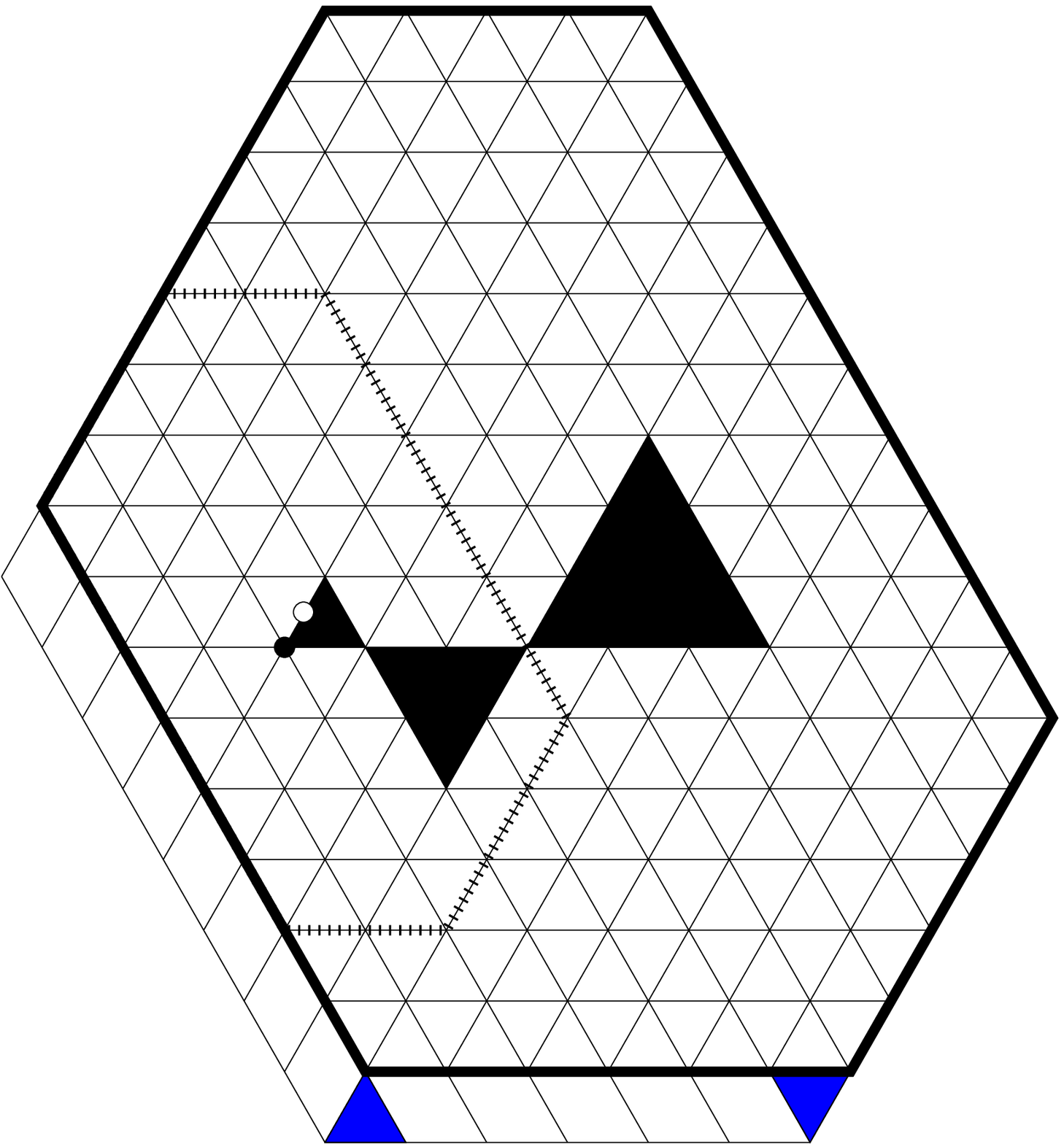}}
\bigskip

\twoline{\mypic{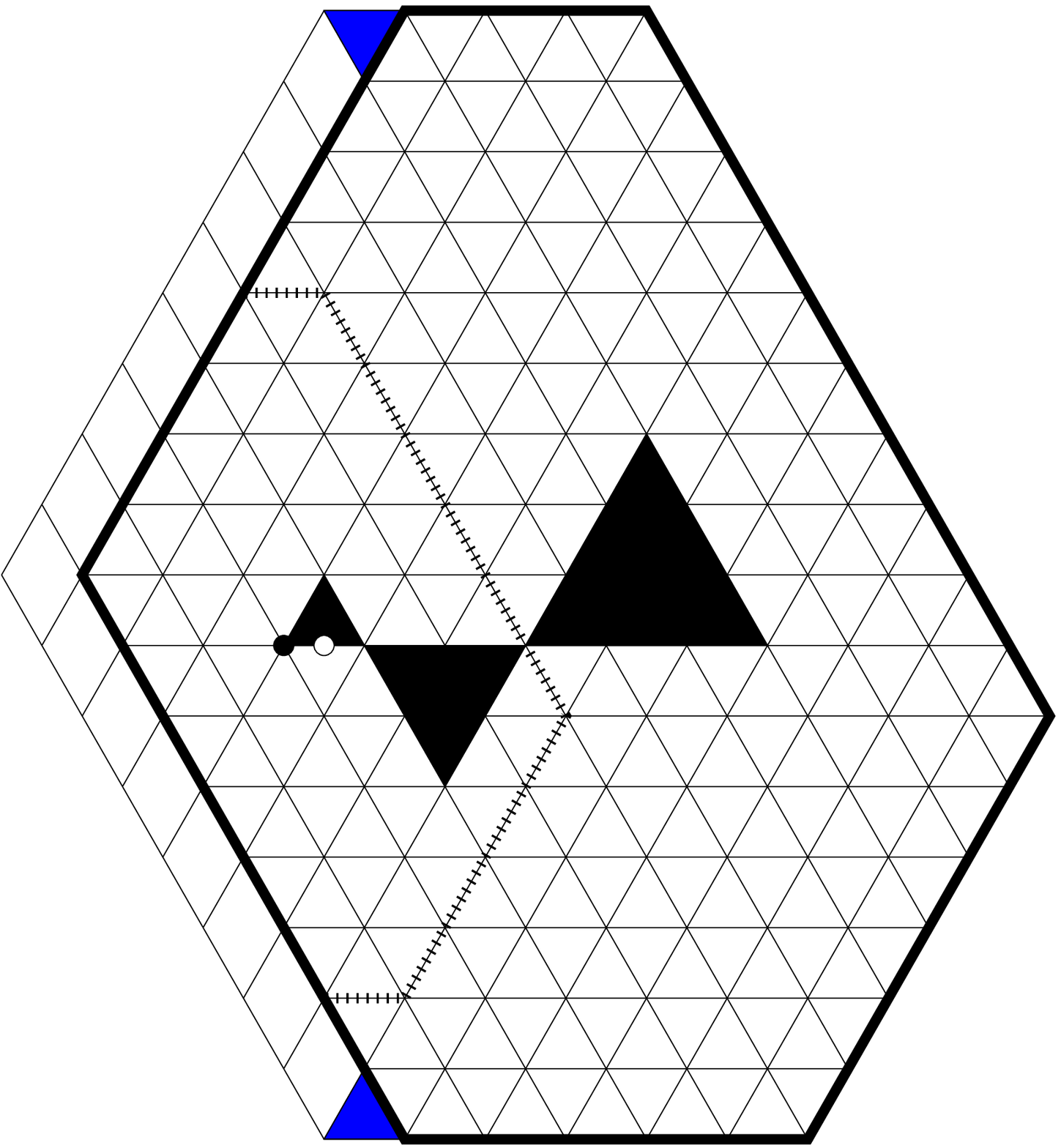}}{\mypic{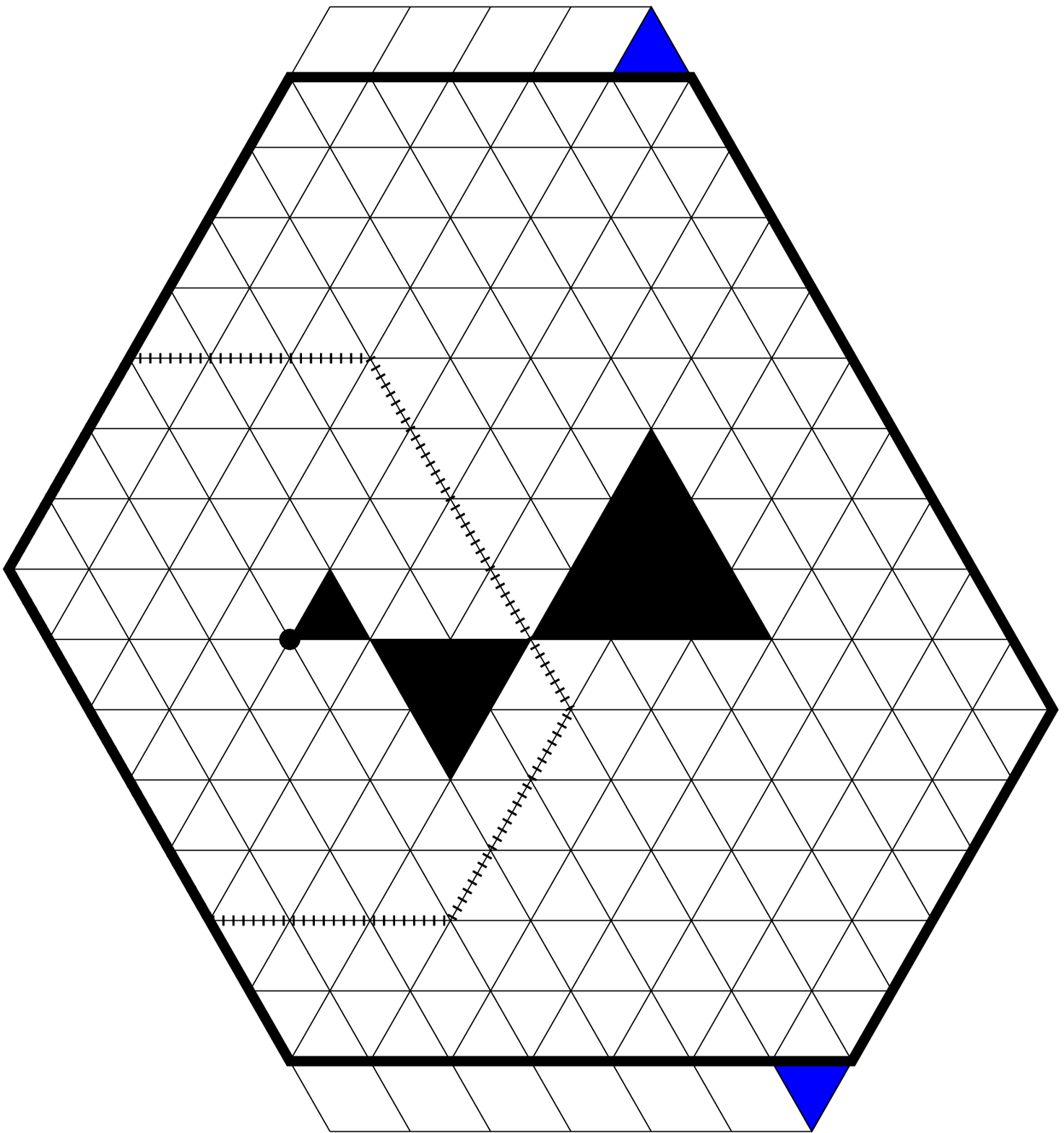}}
\smallpagebreak
\centerline{{\smc Figure {\fca}.}\ The recurrence for the regions $FC^{\bigodot}_{x,y,z}(a_1,\dotsc,a_k)$. Kuo condensation}
\centerline{is applied to the region $FC^{\bigodot}_{2,6,4}(1,2,3)$ (top left) as shown on the top right.}
\endinsert

We recall for convenience the variant of Kuo condensation that we will need in our~proof. 

\proclaim{Theorem {\tca} (Kuo \cite{\Kuo})} Let $G=(V_1,V_2,E)$ be a plane bipartite graph in which $|V_1|=|V_2|$. Let vertices $a$, $b$, $c$ and $d$ appear cyclically on a face of $G$. If $a,c\in V_1$ and $b,d\in V_2$, then
$$
\M(G)\M(G-\{a,b,c,d\})=\M(G-\{a,b\})\M(G-\{c,d\})+\M(G-\{a,d\})\M(G-\{b,c\}).
\tag\eca
$$

\endproclaim

We record here for later use a simple consequence of the definition of the $F$-cored hexagons at the beginning of Section 2.

\proclaim{Lemma {\tcb}} Consider the $F$-cored hexagon $FC_{x,y,z}(a_1,\dotsc,a_k)$, and let $\bar{H}'$ be the translation of its auxilliary hexagon $\bar{H}$ into its eastern corner. Then the rightmost point of the fern core has the same relative position to the center of $\bar{H}'$ as the base of the fern core (i.e., its leftmost point) has to the center of $\bar{H}$ (see the pictures on center right in Figures~{\fcc} and {\fcd} for two illustrations of this). \epf

\endproclaim

{\it Proof of Theorem {\tba}.} We proceed by induction on $x+y+z$. To make it clear what the base cases are, we present first the recurrences that we will use at the induction step. There are four cases, corresponding to the relative parities of $x$, $y$ and $z$.

Suppose $x$, $y$ and $z$ have the same parity, and consider the region $FC^{\bigodot}_{x,y,z}(a_1,\dotsc,a_k)$. Apply Kuo condensation to its planar dual graph, with the vertices $a$, $b$, $c$, $d$ chosen to correspond to the unit triangles on the boundary marked in the top right picture in Figure~{\fca}. The six pictures in Figure {\fca} correspond to the six terms resulting from the identity~(\eca) (rather than showing the dual graphs themselves, Figure {\fca} shows the regions they correspond to). After removing the forced lozenges, we obtain in all cases $F$-cored hexagons, having side-lengths of various relative parities. The auxilliary hexagon (indicated throughout Figure {\fca} by a thick dotted line) helps us read off the precise parameters and kinds of the resulting $F$-cored hexagons. Clearly, the top left picture corresponds to $FC^{\bigodot}_{x,y,z}(a_1,\dotsc,a_k)$ itself. Consider the top right picture. The region left over after removing the forced lozenges is outlined by a thick solid line. Note how in the auxilliary hexagon the $x$-parameter is the same as in the top left picture, but both the $y$- and $z$-parameters are decremented by one unit. Furthermore, the base of the fern is half a unit to the left of the center of the auxilliary hexagon. With this clue in mind, one readily checks that the resulting region is in fact the $F$-cored hexagon $FC^{\leftarrow}_{x,y-1,z-1}(a_1,\dotsc,a_k)$. Similar reasonings lead to the realization that the remaining pictures in Figure {\fca} correspond to the $F$-cored hexagons $FC^{\nwarrow}_{x,y-1,z}(a_1,\dotsc,a_k)$, $FC^{\swarrow}_{x,y,z-1}(a_1,\dotsc,a_k)$, $FC^{\leftarrow}_{x-1,y,z}(a_1,\dotsc,a_k)$, and $FC^{\bigodot}_{x+1,y-1,z-1}(a_1,\dotsc,a_k)$, respectively. Therefore the identity resulting from (\eca) is
$$
\spreadlines{3\jot}
\align
&
\M(FC^{\bigodot}_{x,y,z}(a_1,\dotsc,a_k))\M(FC^{\leftarrow}_{x,y-1,z-1}(a_1,\dotsc,a_k))=
\\
&\ \ \ \ \ \ \ \ \ \ \ \ \ \ \ \ \ \ \ \ \ \ \ \ \ \ \ 
%=
\M(FC^{\nwarrow}_{x,y-1,z}(a_1,\dotsc,a_k))\M(FC^{\swarrow}_{x,y,z-1}(a_1,\dotsc,a_k))
\\
&\ \ \ \ \ \ \ \ \ \ \ \ \ \ \ \ \ \ \ \ \ \ \ \ 
+
\M(FC^{\leftarrow}_{x-1,y,z}(a_1,\dotsc,a_k))\M(FC^{\bigodot}_{x+1,y-1,z-1}(a_1,\dotsc,a_k)).
\tag\ecb
\endalign
$$
When dividing through by the second factor on the left hand side, (\ecb) gives an expression for the number of tilings of $FC^{\bigodot}_{x,y,z}(a_1,\dotsc,a_k)$ in terms of the number of tilings of five other $F$-cored hexagons, all of them having $x$-, $y$- and $z$-parameters that add up to something strictly less than $x+y+z$. This is a recurrence for $\M(FC^{\bigodot}_{x,y,z}(a_1,\dotsc,a_k))$ that we will use at the induction step. Note that in order for all the regions involved in (\ecb) to be defined, we need $x,y,z\geq1$. This shows that the cases $x=0$, $y=0$ and $z=0$ will be base cases for our induction.

\topinsert

\twoline{\mypic{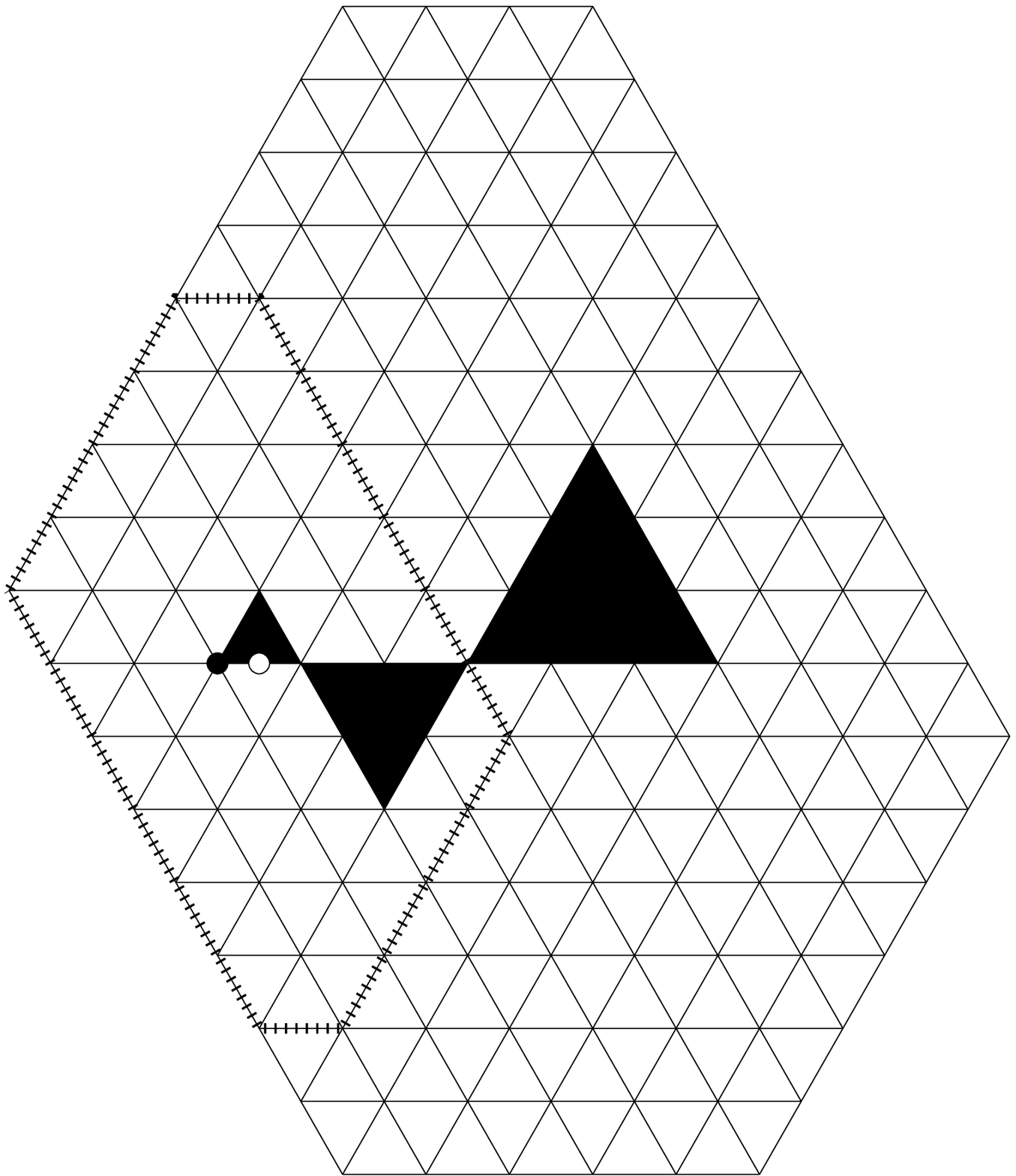}}{\mypic{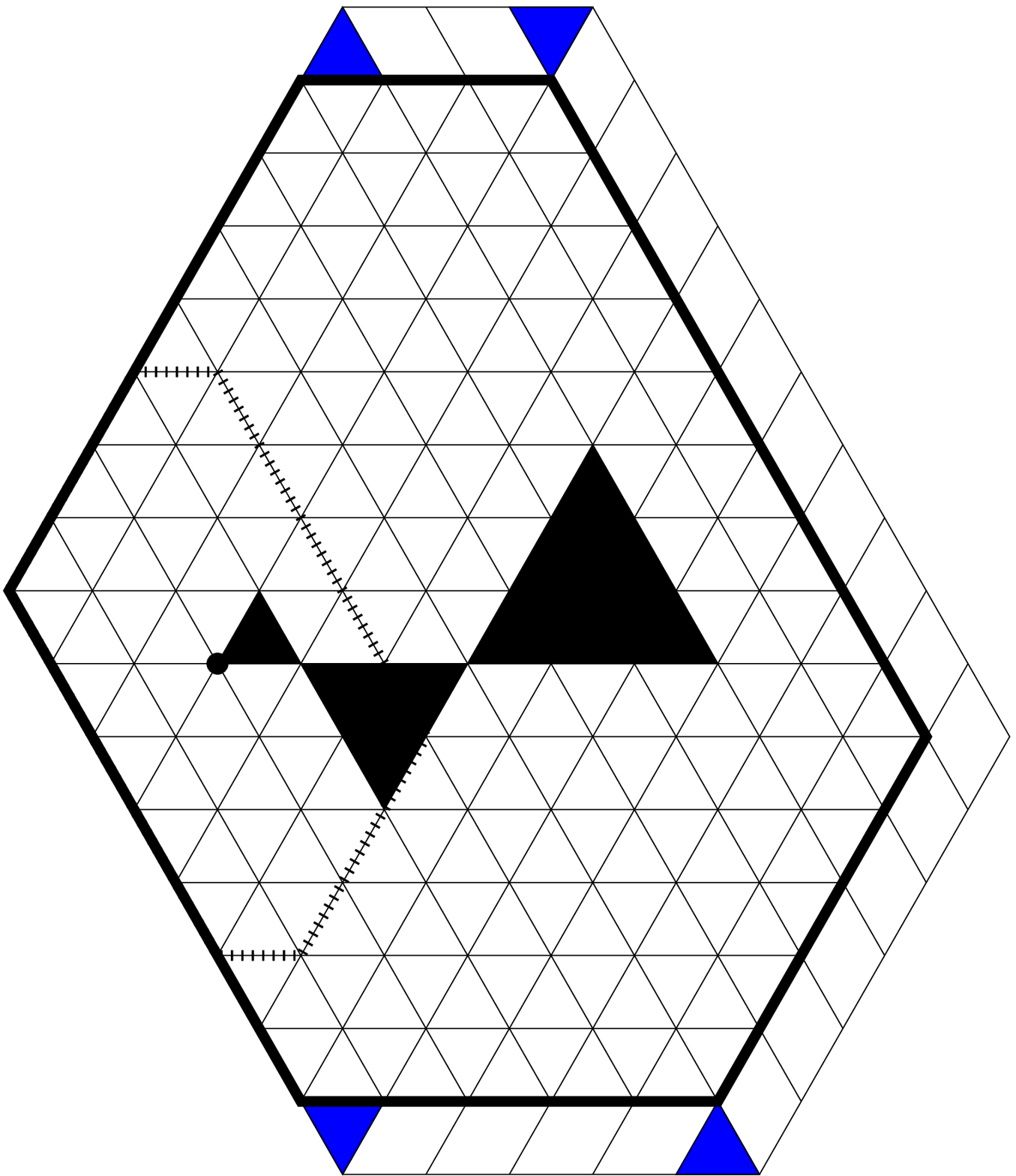}}
\bigskip

\twoline{\mypic{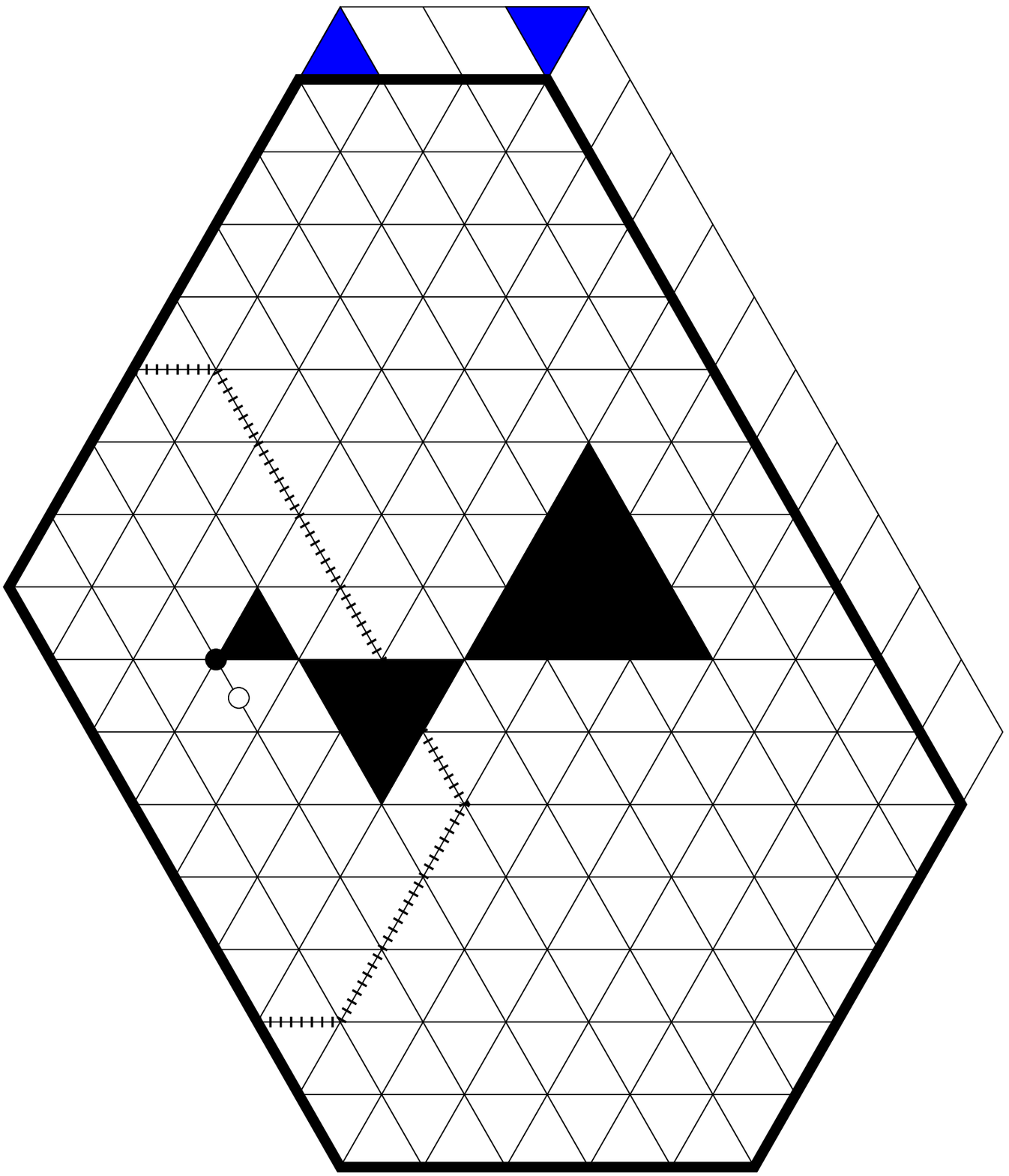}}{\mypic{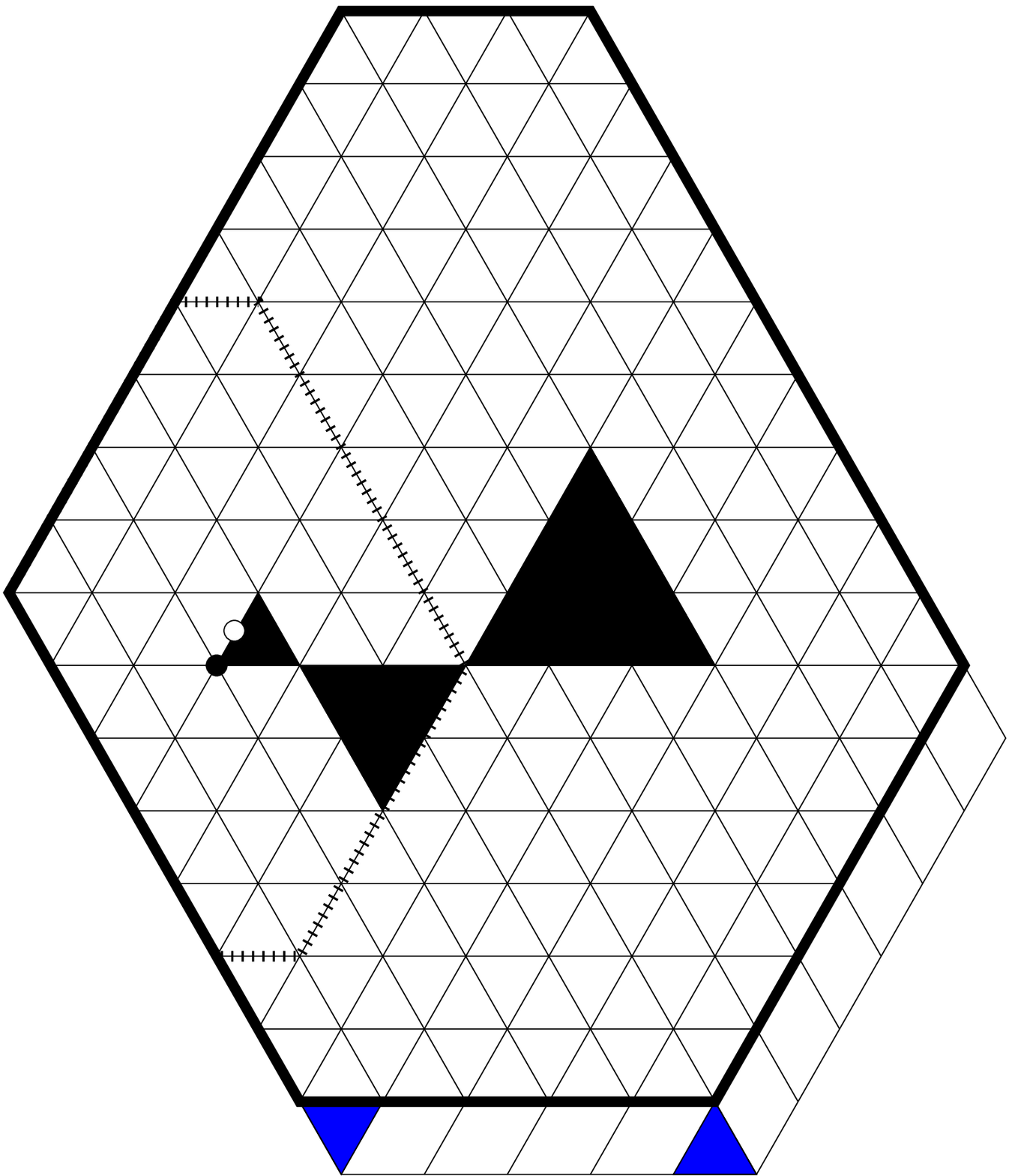}}
\bigskip

\twoline{\mypic{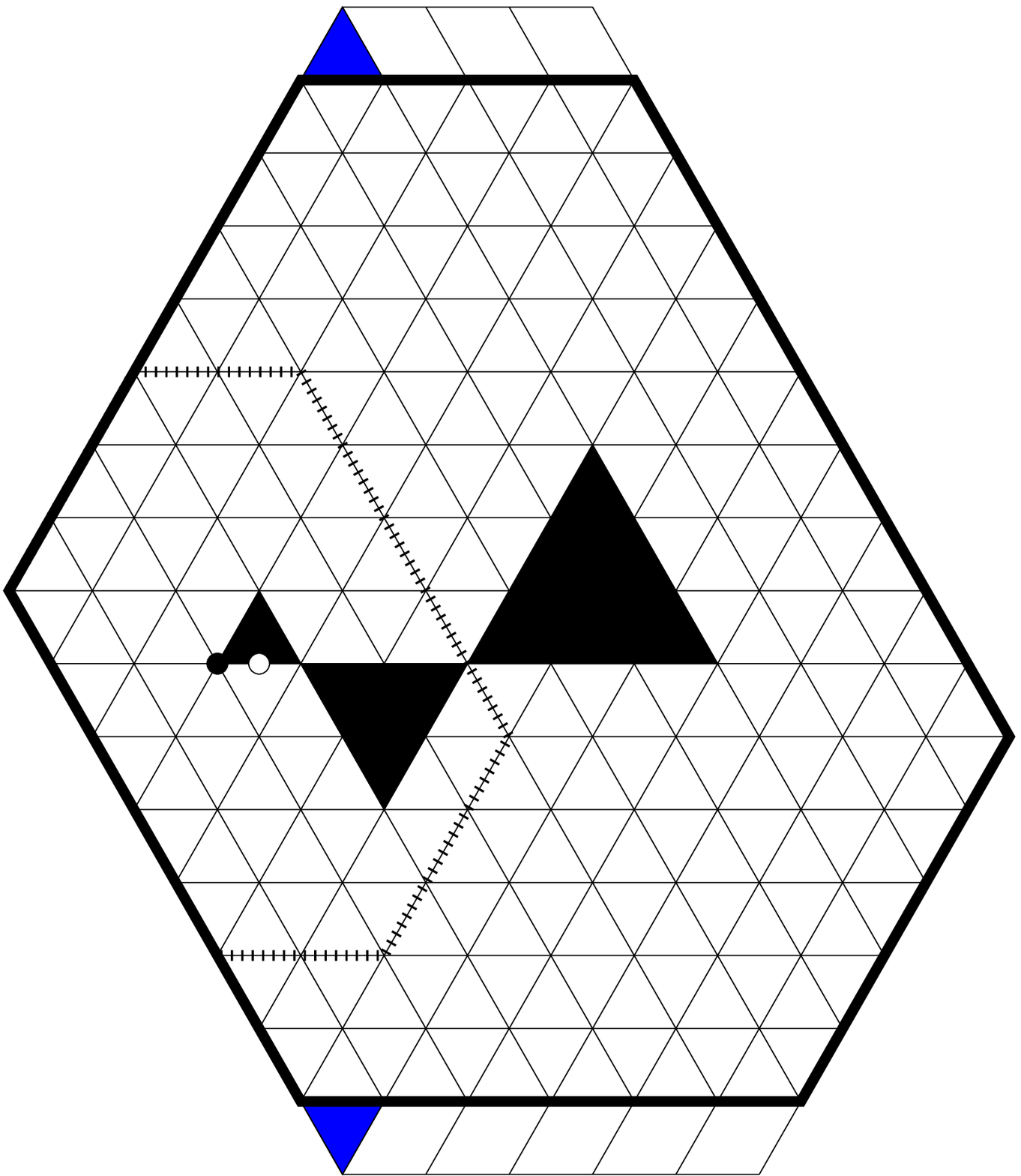}}{\mypic{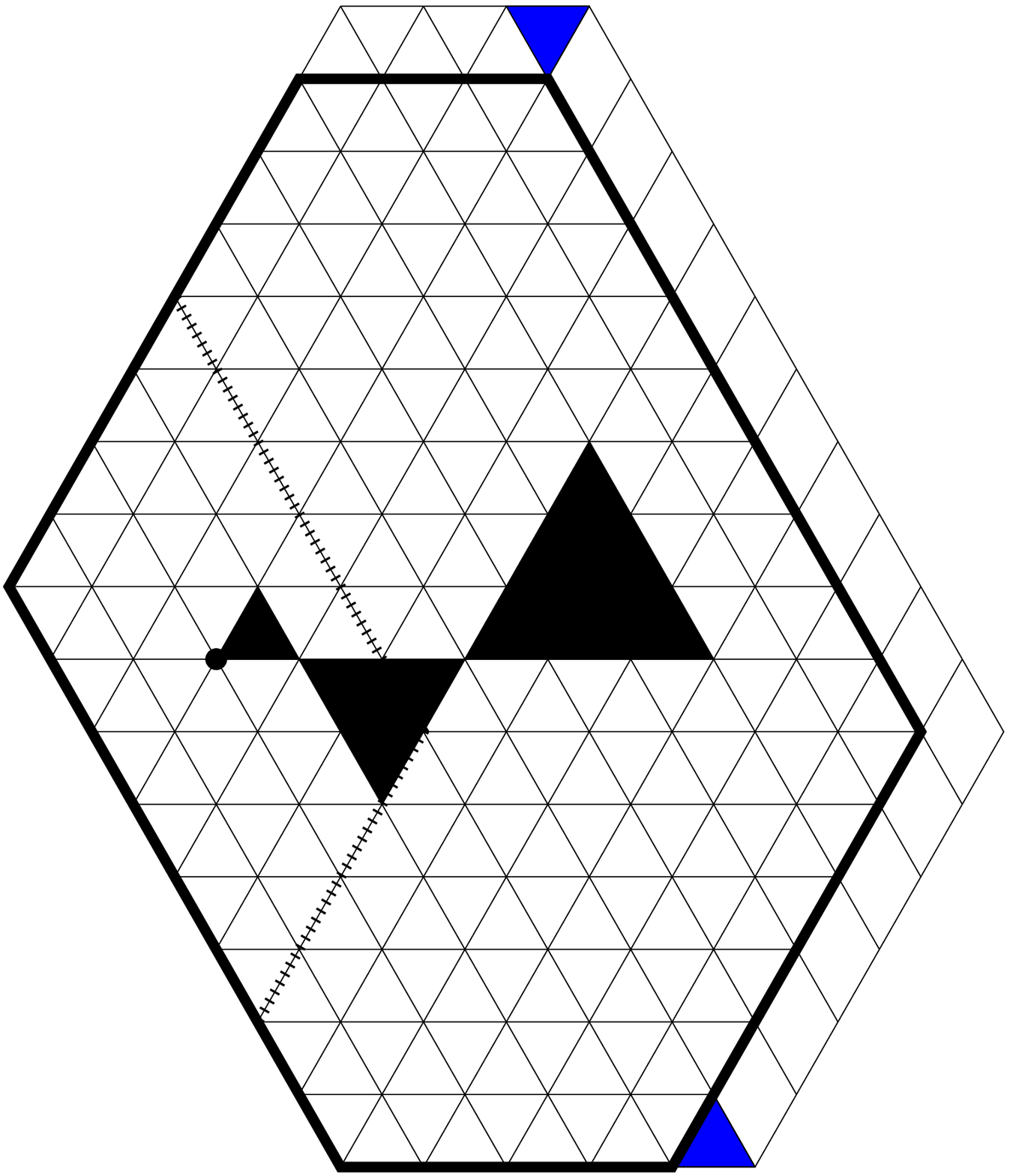}}
\smallpagebreak
\centerline{{\smc Figure {\fcb}.}\ The recurrence for the regions $FC^{\leftarrow}_{x,y,z}(a_1,\dotsc,a_k)$. Kuo condensation}
\centerline{is applied to the region $FC^{\leftarrow}_{1,6,4}(1,2,3)$ (top left) as shown on the top right.}
\endinsert

Since (\ecb) involves also $F$-cored hexagons of types $FC^{\leftarrow}$, $FC^{\swarrow}$ and $FC^{\nwarrow}$, we need recurrences for these types as well.

Consider thus the case when $x$ has parity opposite to $y$ and $z$. Apply Kuo condensation to the dual graph of $FC^{\leftarrow}_{x,y,z}(a_1,\dotsc,a_k)$ as indicated in Figure {\fcb}. Using the same reasoning that led to (\ecb), we obtain
$$
\spreadlines{3\jot}
\align
&
\M(FC^{\leftarrow}_{x,y,z}(a_1,\dotsc,a_k))\M(FC^{\bigodot}_{x,y-1,z-1}(a_1,\dotsc,a_k))=
\\
&\ \ \ \ \ \ \ \ \ \ \ \ \ \ \ \ \ \ \ \ \ \ \ \ \ \ \ 
%=
\M(FC^{\nwarrow}_{x,y,z-1}(a_1,\dotsc,a_k))\M(FC^{\swarrow}_{x,y-1,z}(a_1,\dotsc,a_k))
\\
&\ \ \ \ \ \ \ \ \ \ \ \ \ \ \ \ \ \ \ \ \ \ \ \ 
+
\M(FC^{\leftarrow}_{x+1,y-1,z-1}(a_1,\dotsc,a_k))\M(FC^{\bigodot}_{x-1,y,z}(a_1,\dotsc,a_k)).
\tag\ecc
\endalign
$$
We point out that (\ecb) was obtained by applying Kuo condensation ``at the western corner'' of the $F$-cored hexagon (see Figure {\fca}), i.e. by placing the four unit triangles on the boundary so that the strips of lozenges they force tend to be in the western corner. By contrast, (\ecc) was obtained by applying Kuo condensation at the eastern corner (see Figure {\fcb}). We needed this in the latter case in order to restore as much as possible the central position of the fern core, as in the starting region $FC^{\leftarrow}_{x,y,z}(a_1,\dotsc,a_k)$ the fern core's base was already half a unit to the west of the center of the auxilliary hexagon.

The above two cases had the convenient feature that all the regions that resulted by applying Kuo condensation turned out to be $F$-cored hexagons directly, without the need of rotating or reflecting them. Both remaining cases will involve one region for which these symmetries will be needed. 

Suppose that $z$ has parity opposite to the parity of $x$ and $y$, and consider the region $FC^{\swarrow}_{x,y,z}(a_1,\dotsc,a_k)$. Since the base of the fern core is shifted half a unit in the southwestern lattice direction compared to the center of the auxilliary hexagon, we apply Kuo condensation at the northeastern corner of $FC^{\swarrow}_{x,y,z}(a_1,\dotsc,a_k)$ (in order to bring back the fern core near the center of the resulting regions). Figure {\fcc} shows the resulting regions.

\topinsert

\twoline{\mypic{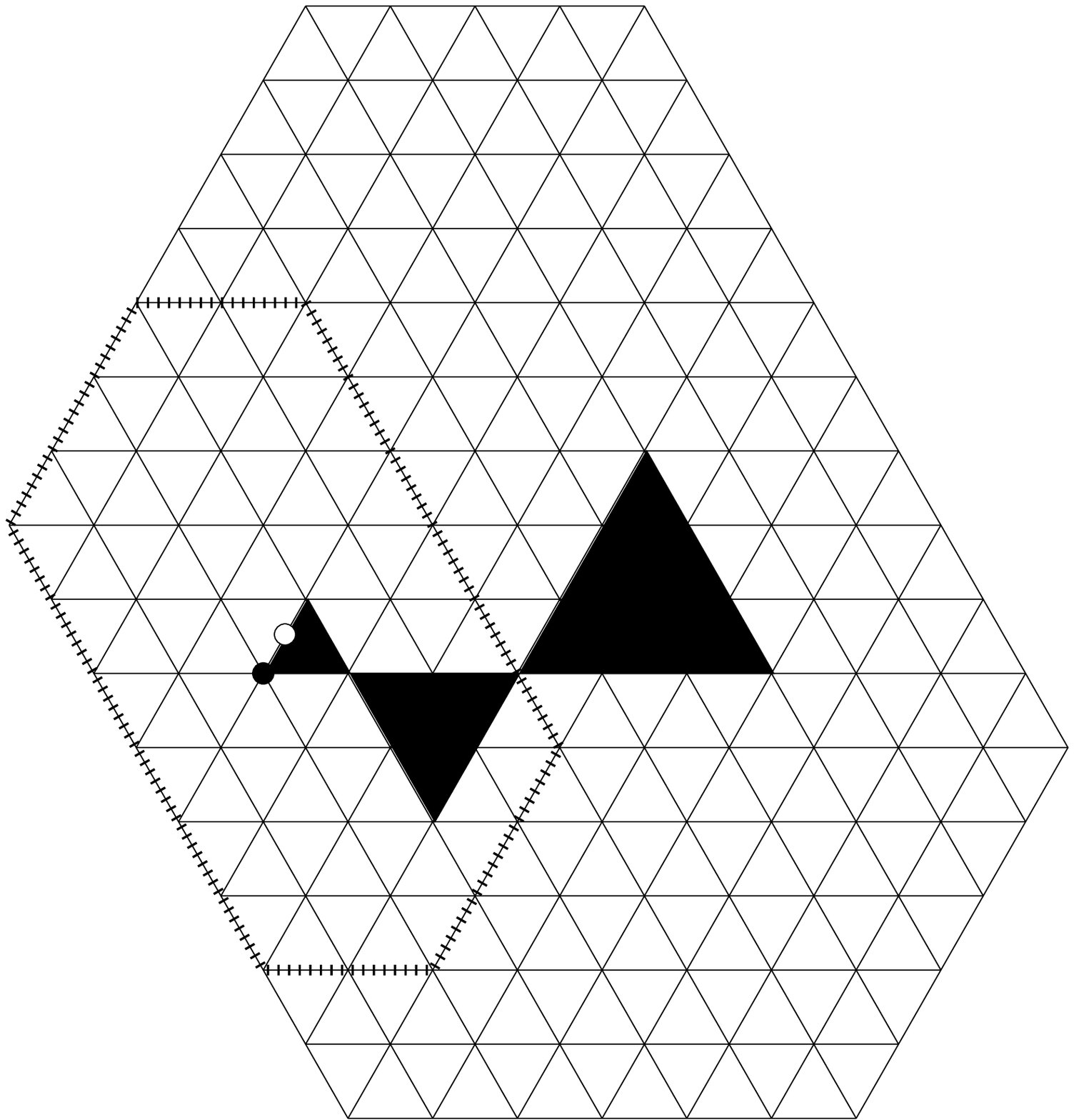}}{\mypic{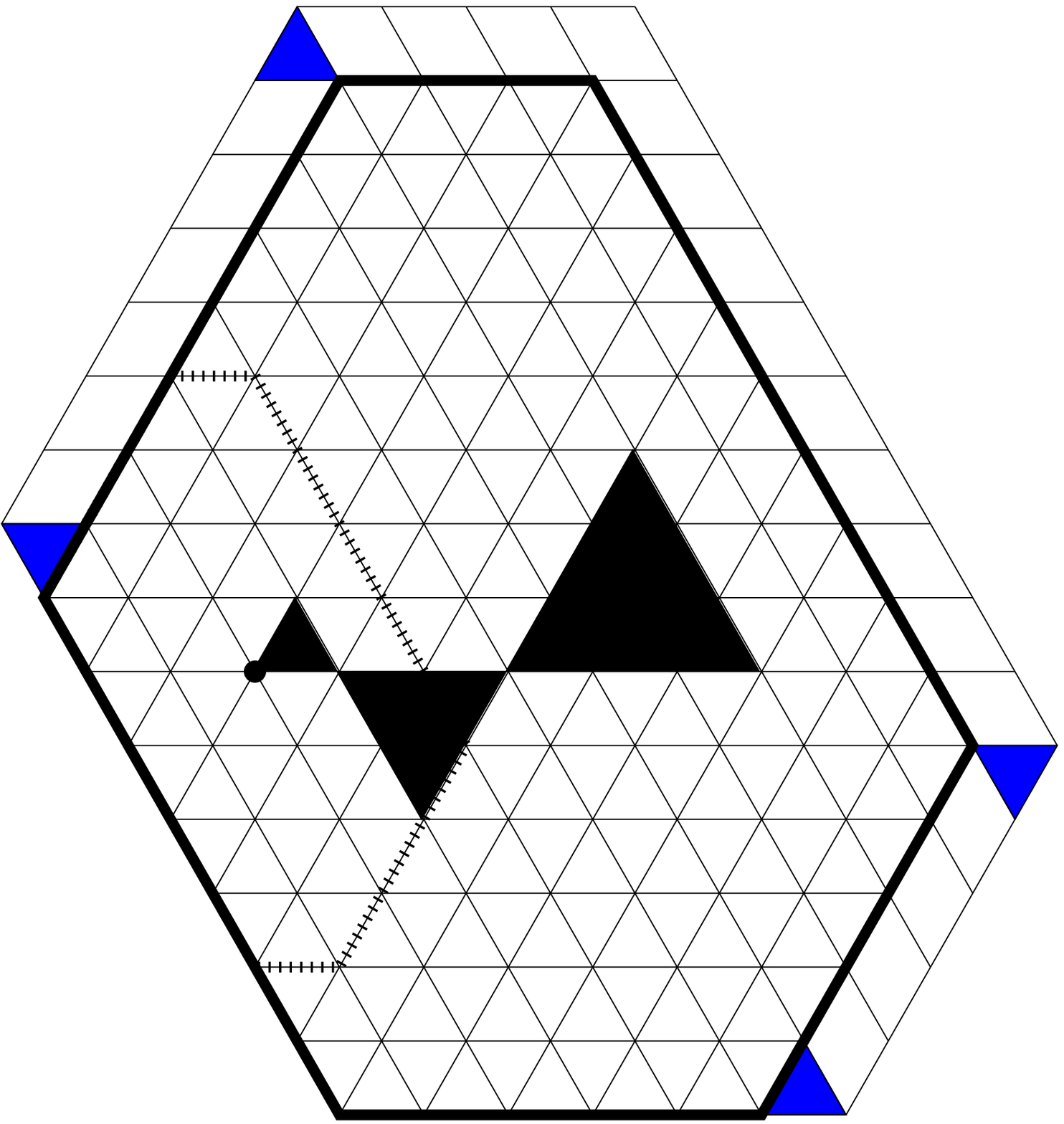}}
\bigskip

\twoline{\mypic{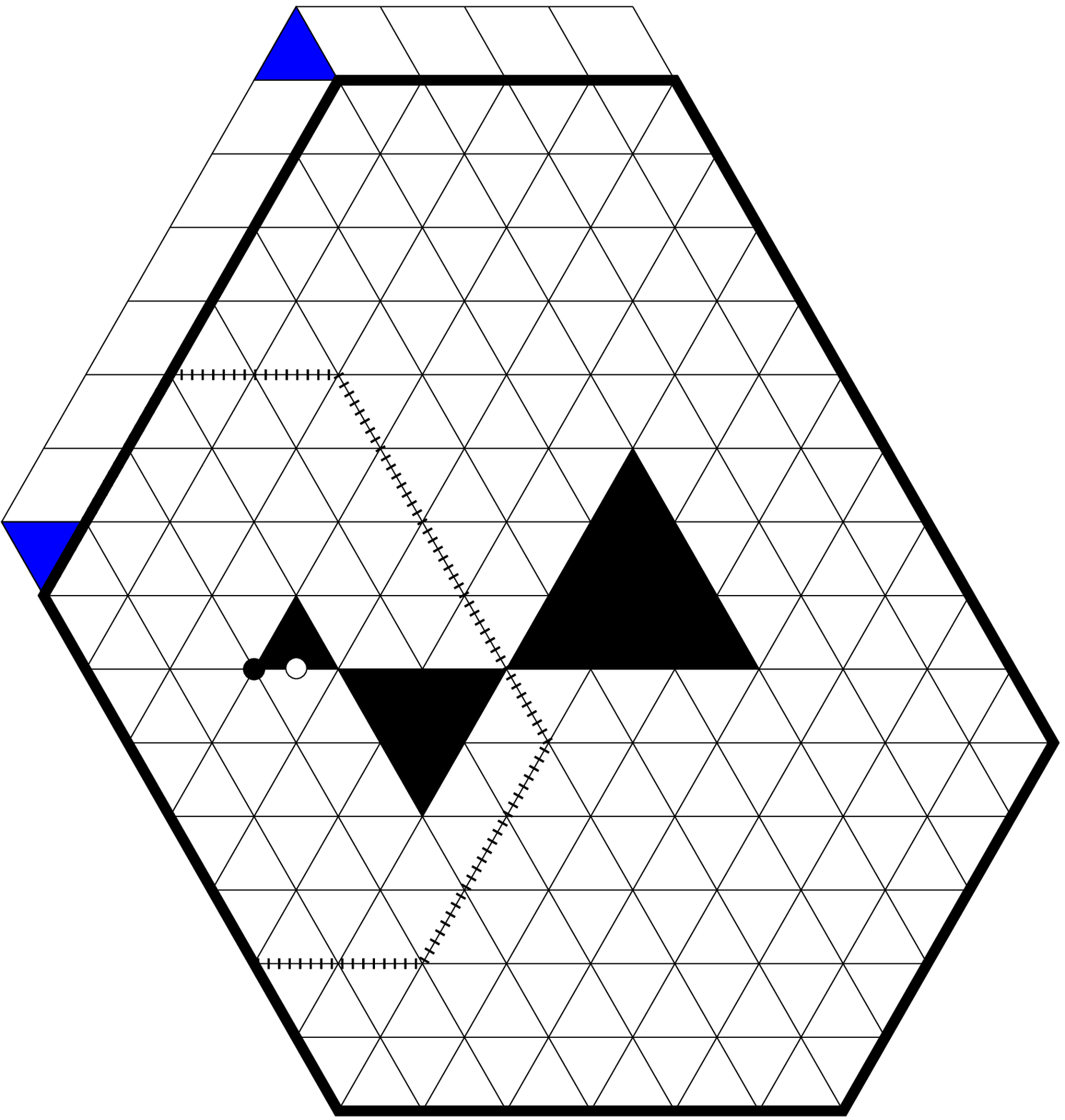}}{\mypic{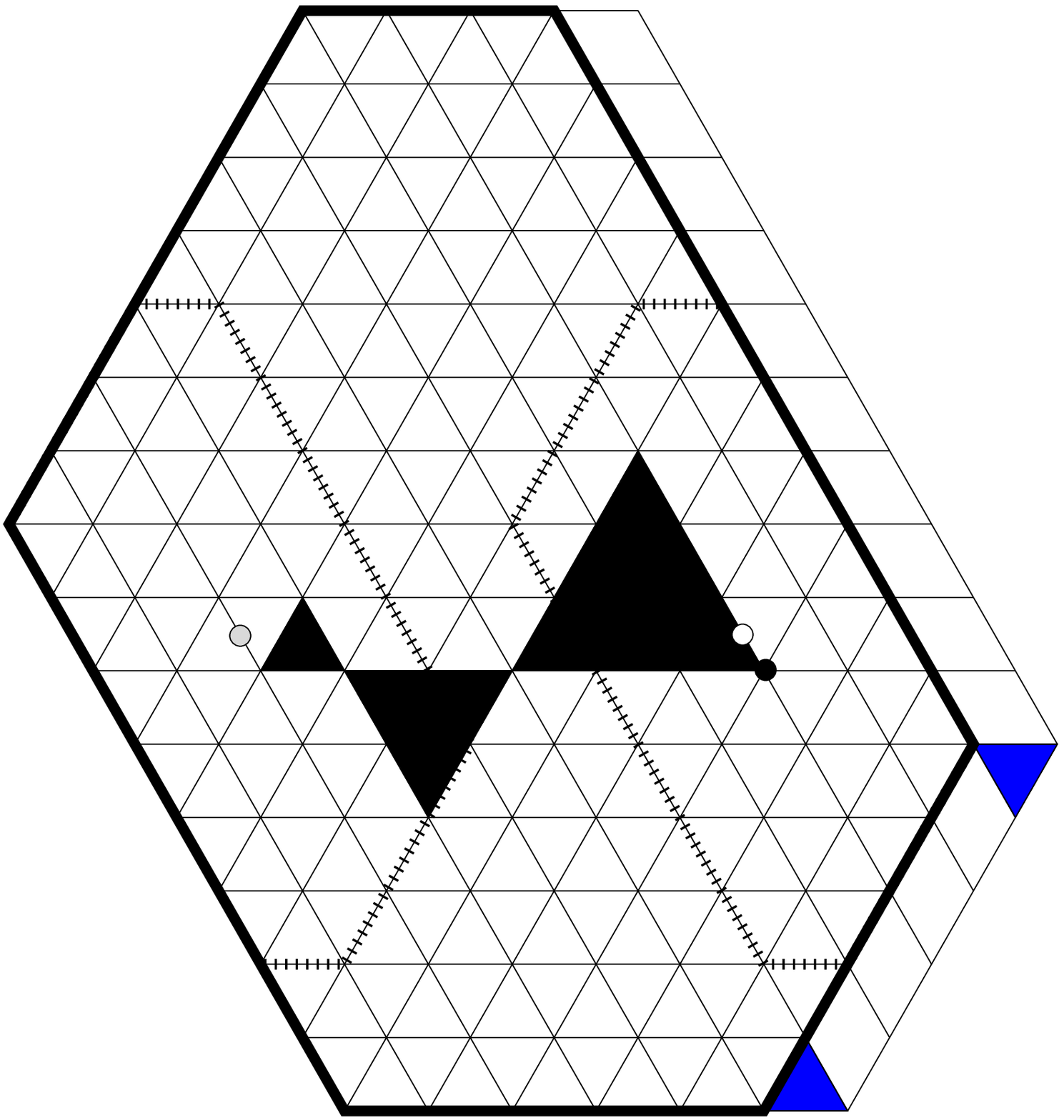}}
\bigskip

\twoline{\mypic{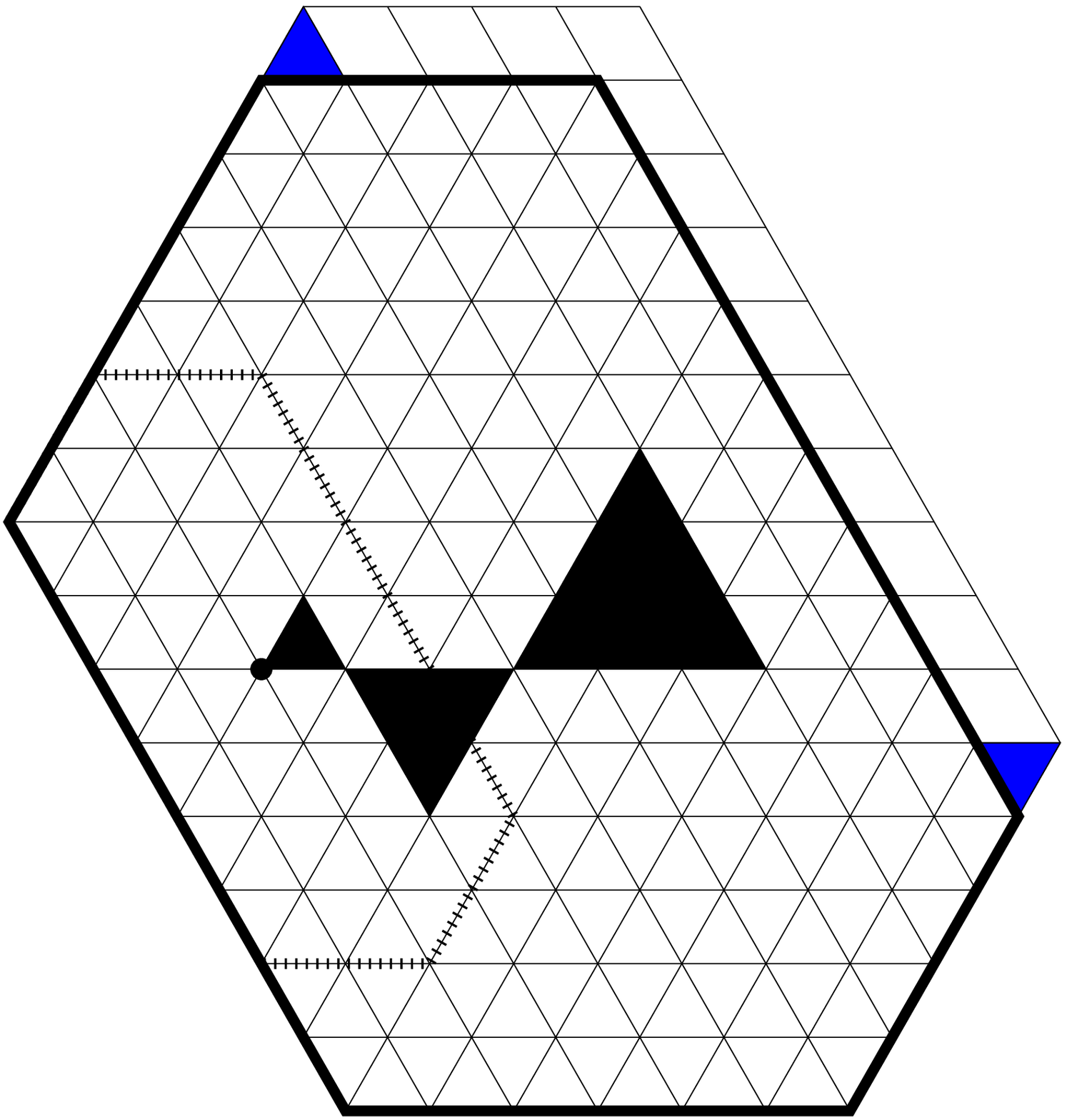}}{\mypic{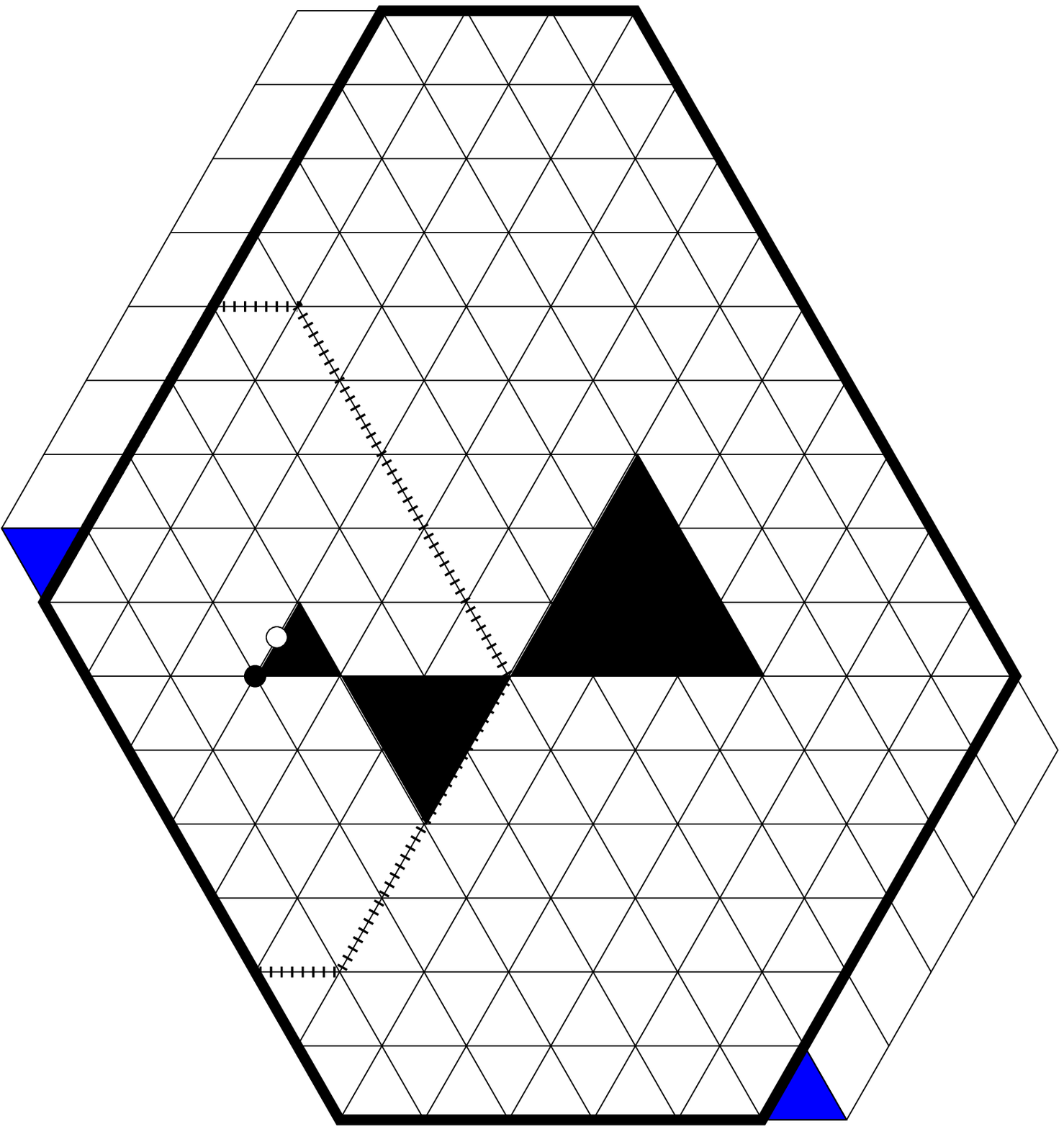}}
\smallpagebreak
\centerline{{\smc Figure {\fcc}.}\ The recurrence for the regions $FC^{\swarrow}_{x,y,z}(a_1,\dotsc,a_k)$. Kuo condensation}
\centerline{is applied to the region $FC^{\swarrow}_{2,6,3}(1,2,3)$ (top left) as shown on the top right.}
\endinsert

With the exception of the center right region in Figure {\fcc}, the resulting regions are readily recognized as $F$-cored hexagons. This can be seen from Figure {\fcc} by looking at the position of the base point of the fern (indicated by a black dot) and the center of the auxilliary hexagon (the auxilliary hexagon is indicated by a thick dotted line, and its center by a white dot): the second, third, fifth and sixth resulting regions are fern-cored hexagons of types $FC^{\bigodot}$, $FC^{\leftarrow}$, $FC^{\bigodot}$ and $FC^{\swarrow}$, respectively, but in the fourth the leftmost point of the fern lies southeast of the center of the auxilliary hexagon, a displacement direction not considered in our family of regions.

The way we resolve this is to view the fourth region after rotating it by $180^\circ$. Once this rotation is performed, by Lemma {\tcb}, the leftmost point of the fern (indicated by a black dot in Figure {\fcc} {\it before} the rotation) is half a unit to the northwest of the center of the auxilliary hexagon (in Figure {\fcc} center right, these are shown before the rotation, hence their location near the eastern corner). As we will see below, the resulting region turns out to be an $F$-cored hexagon of type $FC^{\nwarrow}$. 

There arises, however, a difference between the case when the fern has an even or an odd number of lobes. Namely, if the fern has an even number of lobes, then
after the $180^\circ$ rotation we obtain the $F$-cored hexagon $FC^{\nwarrow}_{x-1,y,z}(a_k,\dotsc,a_1)$ (note the reversal of the order of the lobes). 

On the other hand, if the fern has an odd number of lobes (as in Figure {\fcc}), then after the $180^\circ$ rotation, the fern's lobes start, from left to right, with the first lobe pointing {\it down}, while in the definition of our $F$-cored regions it points up. So we are led to a region $\overline{FC}^{\nwarrow}_{x-1,y,z}(a_k,\dotsc,a_1)$ that we have not considered yet, which differs from $FC^{\nwarrow}_{x-1,y,z}(a_k,\dotsc,a_1)$ in that the lobes point down, up, down, up, and so on, instead of up, down, up, down, and so on.

This is resolved by noting that reflecting the region $\overline{FC}^{\nwarrow}_{x-1,y,z}(a_k,\dotsc,a_1)$ across the horizontal, it becomes our region $FC^{\swarrow}_{x-1,z,y}(a_k,\dotsc,a_1)$.

\smallpagebreak
We obtain therefore the following recurrences:
$$
\spreadlines{3\jot}
\align
&
\M(FC^{\swarrow}_{x,y,z}(a_1,\dotsc,a_k))\M(FC^{\bigodot}_{x-1,y-1,z}(a_1,\dotsc,a_k))=
\\
&\ \ \ \ \ \ \ \ \ \ \ \ \ \ \ \ \ \ \ \ \ \ \ \ \ \ \ 
%=
\M(FC^{\leftarrow}_{x,y-1,z}(a_1,\dotsc,a_k))\M(FC^{\nwarrow}_{x-1,y,z}(a_k,\dotsc,a_1))
%\M(FC^{\swarrow}_{x-1,y,z}\overleftarrow{(a_1,\dotsc,a_k)})
\\
&\ \ \ \ \ \ \ \ \ \ \ \ \ \ \ \ \ \ \ \ \ \ \ \ 
+
\M(FC^{\bigodot}_{x,y,z-1}(a_1,\dotsc,a_k))\M(FC^{\swarrow}_{x-1,y-1,z+1}(a_1,\dotsc,a_k)),
\tag\ecd
\endalign
$$
if the number $k$ of lobes is even, and
$$
\spreadlines{3\jot}
\align
&
\M(FC^{\swarrow}_{x,y,z}(a_1,\dotsc,a_k))\M(FC^{\bigodot}_{x-1,y-1,z}(a_1,\dotsc,a_k))=
\\
&\ \ \ \ \ \ \ \ \ \ \ \ \ \ \ \ \ \ \ \ \ \ \ \ \ \ \ 
%=
\M(FC^{\leftarrow}_{x,y-1,z}(a_1,\dotsc,a_k))\M(FC^{\swarrow}_{x-1,z,y}(a_k,\dotsc,a_1))
\\
&\ \ \ \ \ \ \ \ \ \ \ \ \ \ \ \ \ \ \ \ \ \ \ \ 
+
\M(FC^{\bigodot}_{x,y,z-1}(a_1,\dotsc,a_k))\M(FC^{\swarrow}_{x-1,y-1,z+1}(a_1,\dotsc,a_k)).
\tag\ece
\endalign
$$
if the number of lobes is odd (the only difference between (\ecd) and (\ece) is at the indices of the second regions on their middle lines).

The remaining case is when $y$ has parity opposite to $x$ and $z$. Since now the base of the fern core is to the northwest of the center of the auxilliary hexagon, we apply Kuo condensation at the southeastern corner of the $F$-cored hexagon $FC^{\nwarrow}_{x,y,z}(a_1,\dotsc,a_k)$. Figure {\fcd} shows the resulting regions. A similar reasoning to the one that gave (\ecd) and (\ece) (see Figure {\fcd}) leads to the recurrences
$$
\spreadlines{3\jot}
\align
&
\M(FC^{\nwarrow}_{x,y,z}(a_1,\dotsc,a_k))\M(FC^{\bigodot}_{x-1,y,z-1}(a_1,\dotsc,a_k))=
\\
&\ \ \ \ \ \ \ \ \ \ \ \ \ \ \ \ \ \ \ \ \ \ \ \ \ \ \ 
%=
\M(FC^{\leftarrow}_{x,y,z-1}(a_1,\dotsc,a_k))\M(FC^{\swarrow}_{x-1,y,z}(a_k,\dotsc,a_1))
%\M(FC^{\swarrow}_{x-1,y,z}\overleftarrow{(a_1,\dotsc,a_k)})
\\
&\ \ \ \ \ \ \ \ \ \ \ \ \ \ \ \ \ \ \ \ \ \ \ \ 
+
\M(FC^{\bigodot}_{x,y-1,z}(a_1,\dotsc,a_k))\M(FC^{\nwarrow}_{x-1,y+1,z-1}(a_1,\dotsc,a_k)),
\tag\ecf
\endalign
$$
if the number $k$ of lobes is even, and

\topinsert

\twoline{\mypic{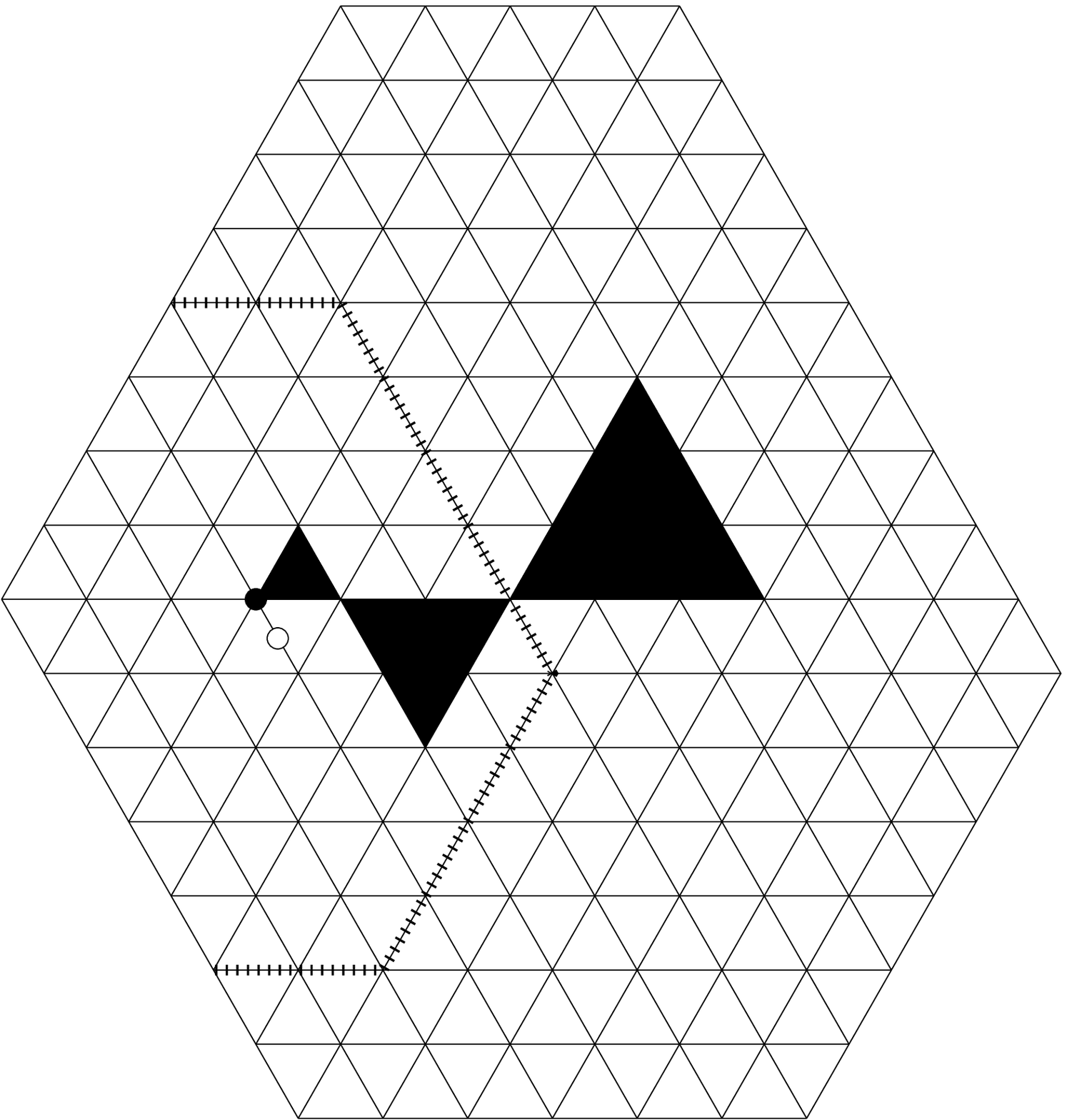}}{\mypic{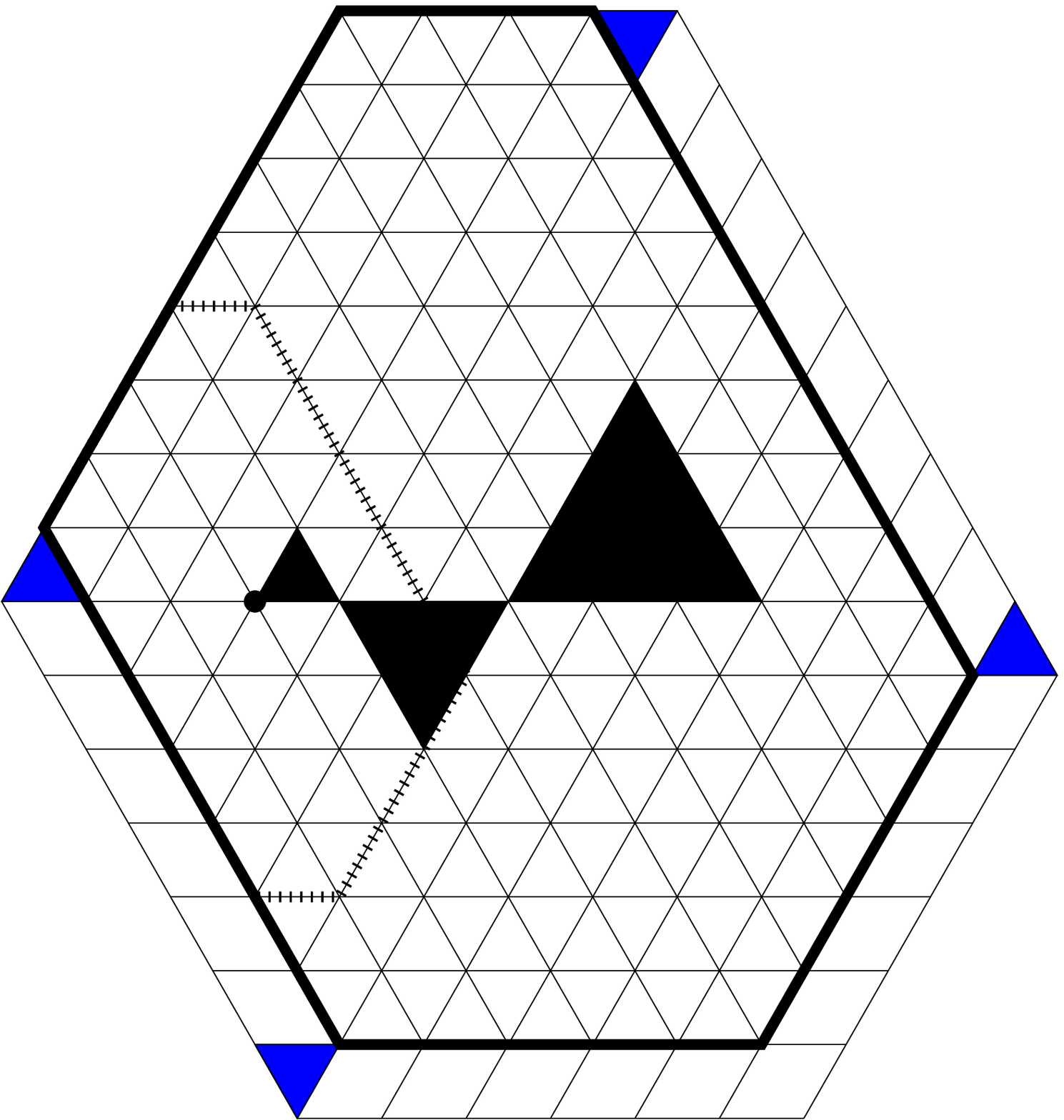}}
\bigskip

\twoline{\mypic{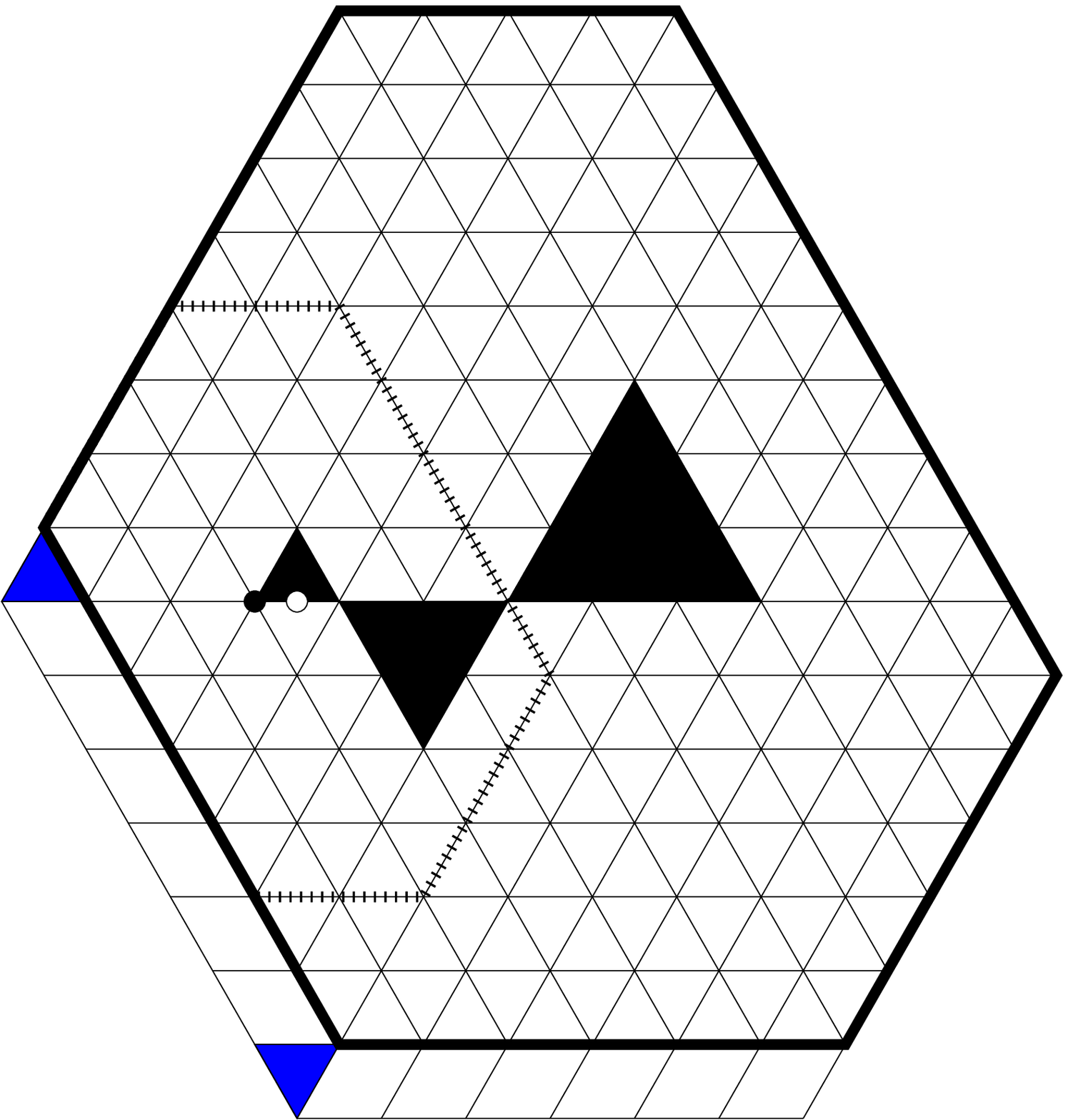}}{\mypic{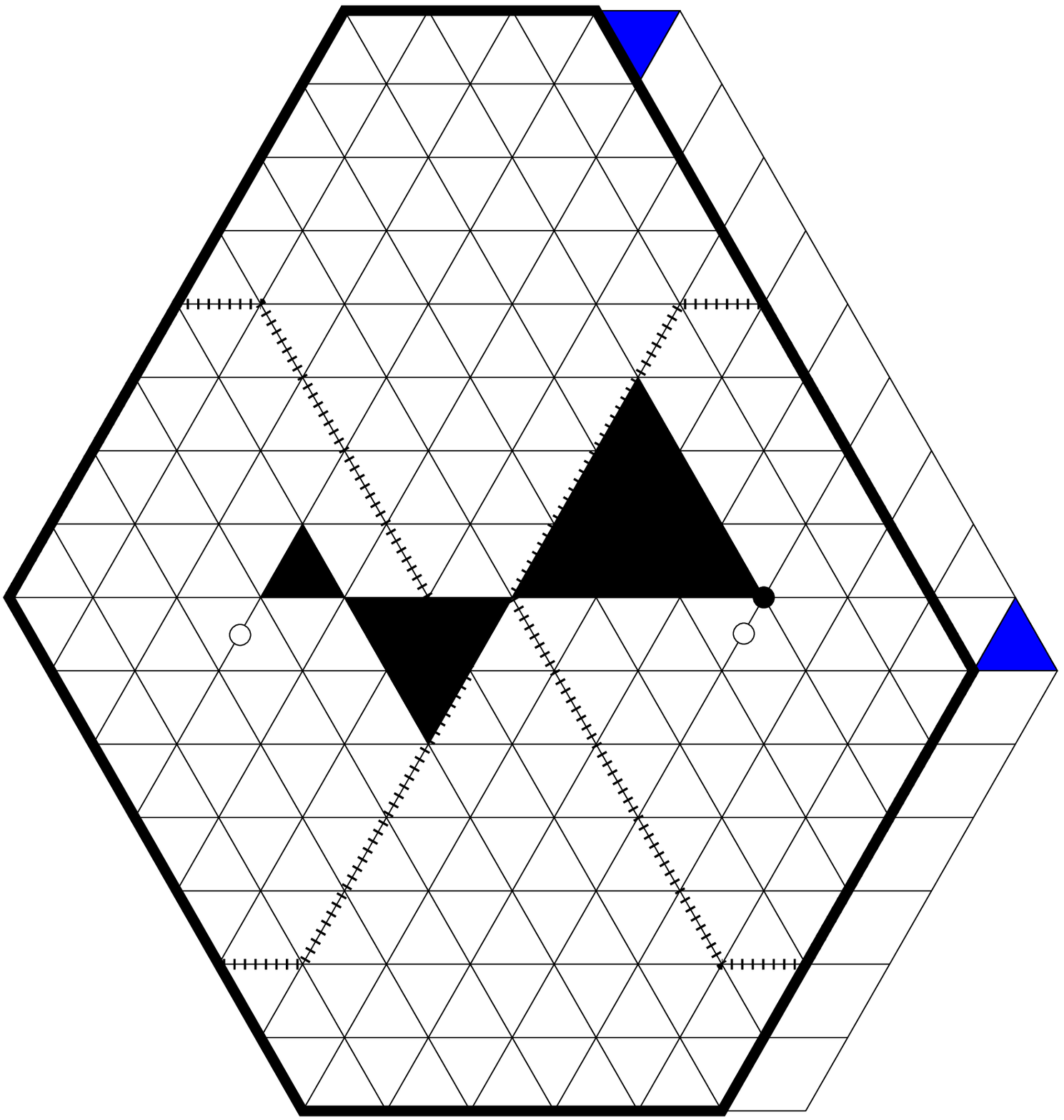}}
\bigskip

\twoline{\mypic{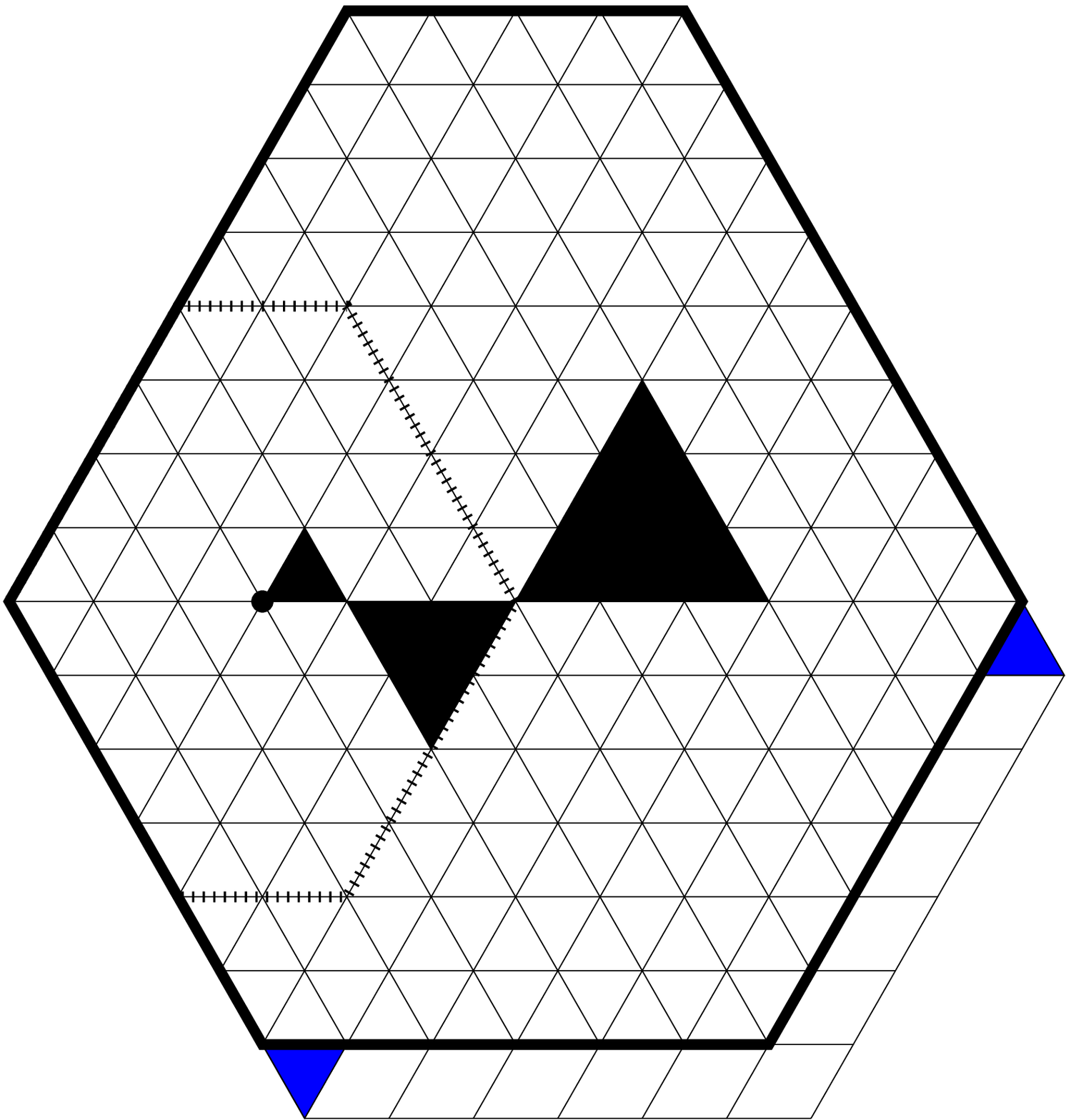}}{\mypic{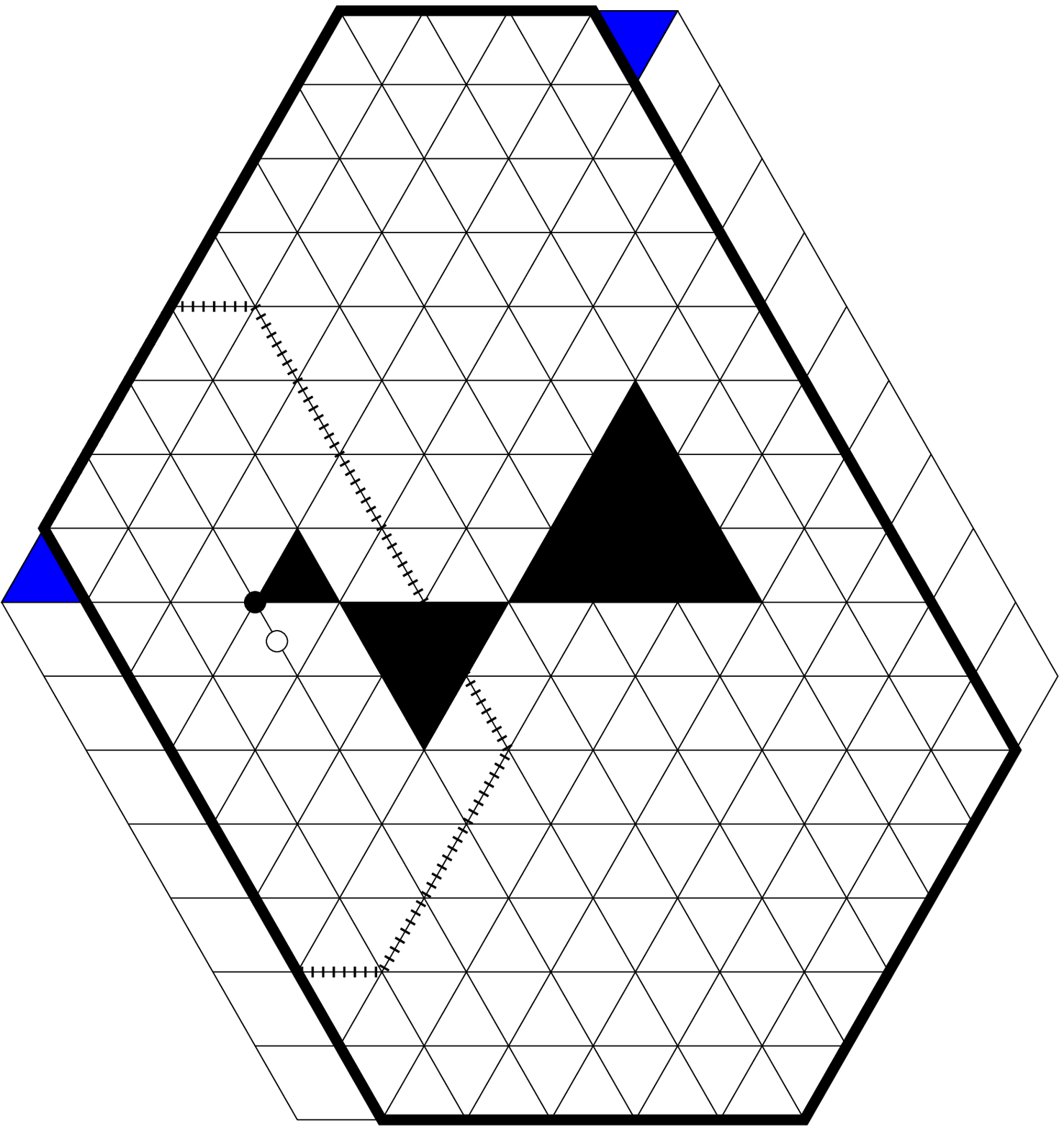}}
\smallpagebreak
\centerline{{\smc Figure {\fcd}.}\ The recurrence for the regions $FC^{\nwarrow}_{x,y,z}(a_1,\dotsc,a_k)$. Kuo condensation}
\centerline{is applied to the region $FC^{\nwarrow}_{2,5,4}(1,2,3)$ (top left) as shown on the top right.}
\endinsert

$$
\spreadlines{2\jot}
\align
&
\M(FC^{\nwarrow}_{x,y,z}(a_1,\dotsc,a_k))\M(FC^{\bigodot}_{x-1,y,z-1}(a_1,\dotsc,a_k))=
\\
&\ \ \ \ \ \ \ \ \ \ \ \ \ \ \ \ \ \ \ \ \ \ \ \ \ \ \ 
%=
\M(FC^{\leftarrow}_{x,y,z-1}(a_1,\dotsc,a_k))\M(FC^{\nwarrow}_{x-1,z,y}(a_k,\dotsc,a_1))
%\M(FC^{\swarrow}_{x-1,y,z}\overleftarrow{(a_1,\dotsc,a_k)})
\\
&\ \ \ \ \ \ \ \ \ \ \ \ \ \ \ \ \ \ \ \ \ \ \ \ 
+
\M(FC^{\bigodot}_{x,y-1,z}(a_1,\dotsc,a_k))\M(FC^{\nwarrow}_{x-1,y+1,z-1}(a_1,\dotsc,a_k)),
\tag\ecg
\endalign
$$
if the number of lobes is odd (again, the only difference between (\ecf) and (\ecg) is at the indices of the second regions on their middle lines).

To summarize, (3.2), (3.3), (3.4) and (3.6) give the desired set of recurrences when the number of lobes of the fern is even, and (3.2), (3.3), (3.5) and (3.7) give the needed recurrences when the number of lobes is odd. These will be used in the induction step.

Note that a common feature of all six recurrences (3.2)--(3.7) is that all the regions they involve are well defined provided $x,y,z\geq1$. Therefore the base cases of our induction on $x+y+z$ will be the three cases when $x=0$, $y=0$, or $z=0$.

\topinsert
\twoline{\mypic{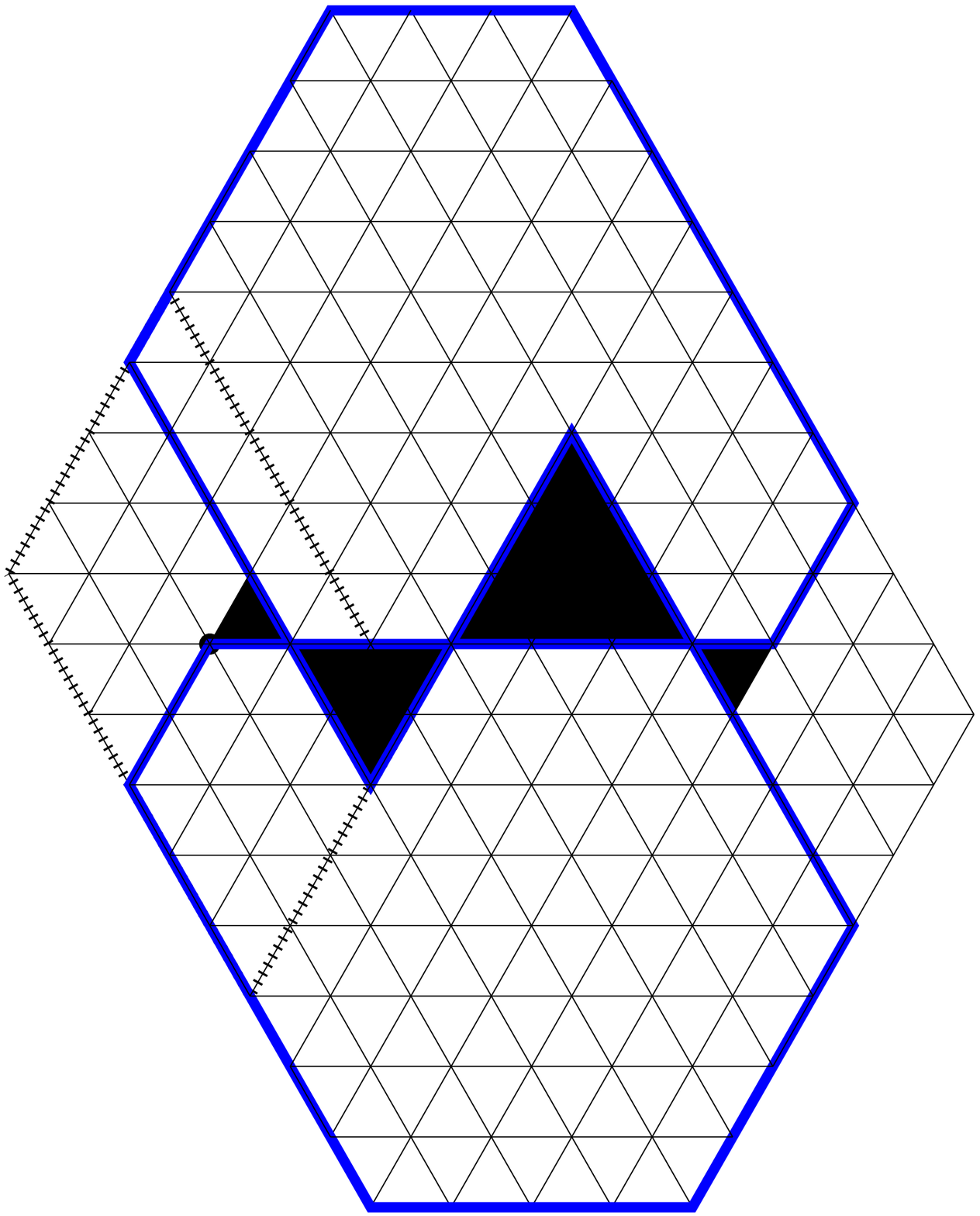}}{\mypic{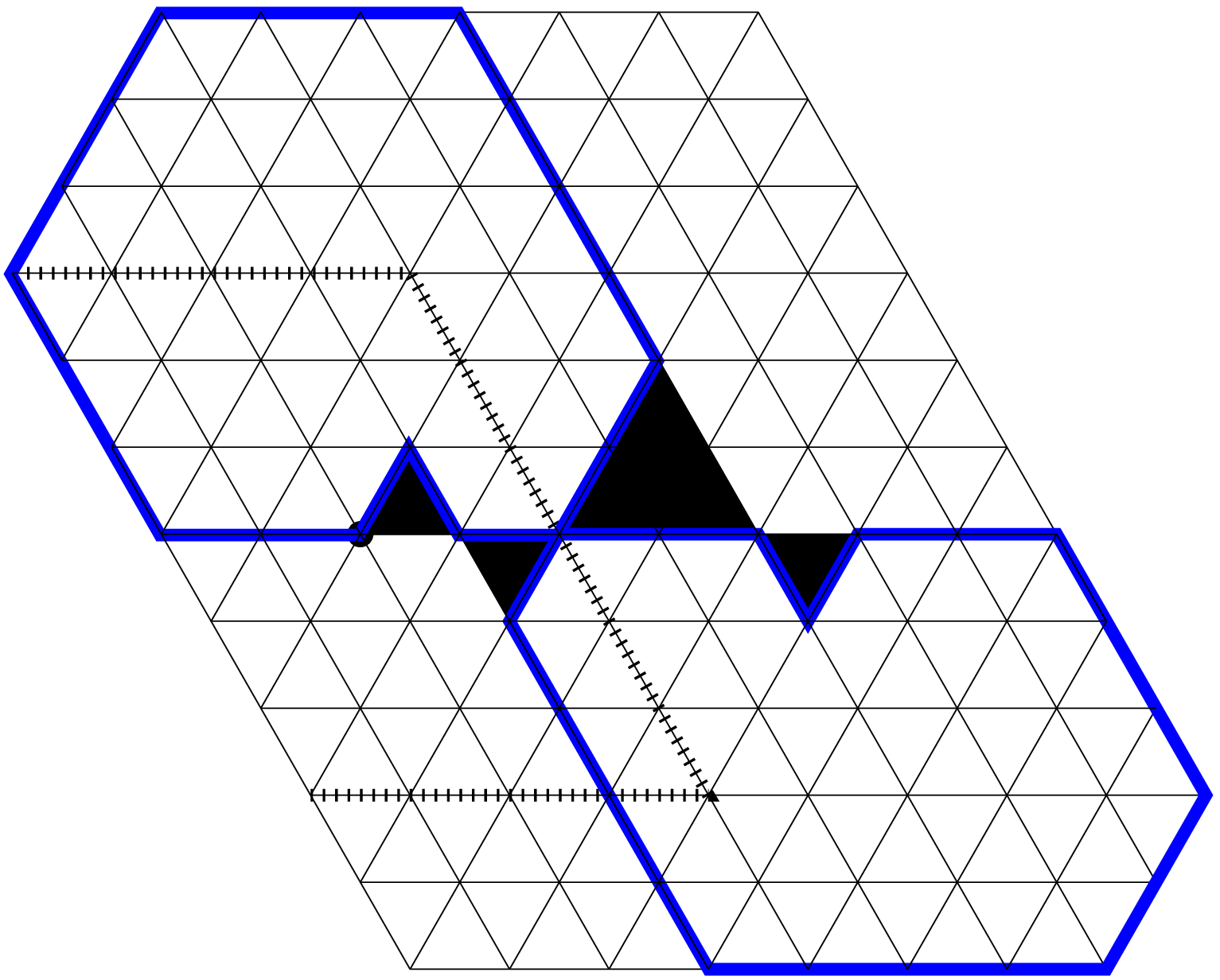}}
\smallpagebreak
\centerline{{\smc Figure {\fce}.}\ The regions $F^{\bigodot}_{0,6,4}(1,2,3,1)$) (left) and $F^{\bigodot}_{4,6,0}(1,1,2,1)$ (right).}
\endinsert

If $z=0$, then the $F$-cored hexagons of type $FC^{\bigodot}$ look as illustrated on the right in Figure {\fce}. As $z=0$, the length of the northwestern side is $a_1+a_3+\cdots$, which matches the total length of the sides of the fern facing northwest. This implies that the upper hexagon $H_1$ outlined in Figure {\fce}, right by a thick contour is internally tiled in each lozenge tiling of $FC^{\bigodot}_{x,y,0}(a_1,\dotsc,a_k)$. A similar argument shows that the lower hexagon $H_2$ outlined in Figure {\fce}, right by a thick contour is also internally tiled. Since on the portion of the $F$-cored hexagon outside these two hexagons the tiling is forced, we obtain that
$$
\M(FC^{\bigodot}_{x,y,0}(a_1,\dotsc,a_k))
=
\M(H_1)\M(H_2).
\tag\eci
$$
However, $H_1$ and $H_2$ are semihexagons of the type covered by formula (\ead). More precisely, when the number of lobes is even, we have 
$$
H_1=S(y/2,x/2,a_1,\dotsc,a_{k-1})\tag\ecj
$$
and 
$$
H_2=S(a_2,\dotsc,a_k,x/2,y/2)
\tag\eck
$$
(note that since $z=0$ and $x,y,z$ have the same parity --- as we are in the $FC^{\bigodot}$ type --- both $x/2$ and $z/2$ are integers).
Plugging in (\ecj) and (\eck) in (\eci) and using formula (\ead), one verifies that the resulting formula for $\M(FC^{\bigodot}_{x,y,0}(a_1,\dotsc,a_k))$ agrees with the one provided by (\ebb). The case of an odd number of lobes, as well as the remaining types of $F$-cored hexagons, are handled in the same way. This also proves the base case $y=0$, by symmetry. The case $x=0$ is analogous.

Assume therefore from now on that $x,y,z\geq1$, and suppose that formula (\ebb) holds for all $F$-cored hexagons for which the sum of the $x$-, $y$- and $z$-parameters is strictly less than $x+y+z$. Consider the $F$-cored hexagon $FC_{x,y,z}(a_1,\dotsc,a_k)$. We remind the reader that its type is determined by the relative parities of $x$, $y$ and $z$, so since $x$, $y$ and $z$ are given, this has one specific type --- either $F^{\bigodot}$, $F^{\leftarrow}$, $F^{\swarrow}$ or $F^{\nwarrow}$. Thus precisely one of the six recurrences (\ecb)--(\ecg) applies to it. Apply it. Then $\M(FC_{x,y,z}(a_1,\dotsc,a_k))$ gets expressed in terms of numbers of tilings of $F$-cored hexagons for which the sum of the $x$-, $y$- and $z$-parameters is strictly less than $x+y+z$ (a common feature to (\ecb)--(\ecg)), for which formula (\ebb) applies by the induction hypothesis. To complete the proof it suffices to check that the resulting expressions agree with what the right hand side of (\ebb) provides. This verification is carried out in Section 4. This completes the proof. \epf

\bigskip

%(?) Lemma. Auxilliary hexagon translated to eastern corner has (1) base at same height as original aux hex and (2) center at same relative position vs rightmost point of fern as center of original aux hex vs leftmost point of fern (i.e. base point)

%\bigskip
%(Basically needed for our claims about fourth pictures in Figures {\fcc} and {\fcd}).

\mysec{4. Verifying that formula (\ebb) satisfies recurrences (\ecb)--(\ecg)}

In this section we verify that the product expression for $\M(FC_{x,y,z}(a_1,\dotsc,a_k))$ yielded by equation (\ebb) --- which we will denote by $f_{x,y,z}(a_1,\dotsc,a_k)$ --- satisfies recurrences (\ecb)--(\ecg). 

By the above definition, we have
$$
f_{x,y,z}(a_1,\dotsc,a_k)
=
\M(FC_{x,y,z}(o,e))
\,g_{x,y,z}(a_1,\dotsc,a_k),
\tag\eda
$$
where 
%$$
%\spreadlines{3\jot}
%\align
%o&:=a_1+a_3+\cdots\tag\edb
%\\
%e&:=a_2+a_4+\cdots\tag\edc
%\endalign
%$$
%and
$$
\spreadlines{3\jot}
\align
&
g_{x,y,z}(a_1,\dotsc,a_k)
=
s(a_1,\dotsc,a_{k-1})s(a_2,\dotsc,a_k)
\\
&\ \ \ \ \ \ \ \ \ \ \ \ \ \ \ \ \ \ \ \ \ \ \ 
\times
\frac
{\h(\lfloor \frac{x+z}{2}\rfloor+a_1+a_3+\cdots)}
{\h(\lfloor \frac{x+y}{2}\rfloor+a_1+a_3+\cdots)}
\frac
{\h(\lceil \frac{x+z}{2}\rceil+a_2+a_4+\cdots)}
{\h(\lceil \frac{x+y}{2}\rceil+a_2+a_4+\cdots)}
\\
&\ \ \ \ \ \ \ \ \ \ 
\times
\prod_{1\leq 2i+1\leq k} 
\frac
{\h(\lfloor \frac{x+y}{2}\rfloor+a_1+\cdots+a_{2i+1})}
{\h(\lfloor \frac{x+z}{2}\rfloor+a_1+\cdots+a_{2i+1})}
\frac
{\h(\lceil \frac{x+y}{2}\rceil+\overline{a_1+\cdots+a_{2i+1}})}
{\h(\lceil \frac{x+z}{2}\rceil+\overline{a_1+\cdots+a_{2i+1}})}
\\
&\ \ \ \ \ \ \ \ \ \ \ 
\times
\prod_{1< 2i< k} 
\frac
{\h(\lfloor \frac{x+z}{2}\rfloor+a_1+\cdots+a_{2i})}
{\h(\lfloor \frac{x+y}{2}\rfloor+a_1+\cdots+a_{2i})}
\frac
{\h(\lceil \frac{x+z}{2}\rceil+\overline{a_1+\cdots+a_{2i}})}
{\h(\lceil \frac{x+y}{2}\rceil+\overline{a_1+\cdots+a_{2i}})},
\tag\edb
\endalign
$$
with $\M(FC_{x,y,z}(a,b))$ given by formulas (\ebg)--(\ebhb), and $s(b_1,\dotsc,b_l)$ by formula~(\ead).

The calculations needed to check that $f_{x,y,z}(a_1,\dotsc,a_k)$ satisfies the six recurrences (\ecb)--(\ecg) are quite similar. We start by presenting in detail these calculations for recurrence (\ecb). For this, we need to check that
$$
\spreadlines{3\jot}
\align
&
f^{\bigodot}_{x,y,z}(a_1,\dotsc,a_k)f^{\nwarrow}_{x,y-1,z-1}(a_1,\dotsc,a_k)=
%\\
%&\ \ \ \ \ \ \ \ \ \ \ \ \ \ \ \ \ \ \ \ \ \ \ \ \ \ \ 
%=
f^{\nwarrow}_{x,y-1,z}(a_1,\dotsc,a_k)f^{\swarrow}_{x,y,z-1}(a_1,\dotsc,a_k)
\\
&\ \ \ \ \ \ \ \ \ \ \ \ \ \ \ \ \ \ \ \ \ \ \ \ \ \ \ \ \ \ \ \ \ \ \ \ \ \ \ 
+
f^{\leftarrow}_{x-1,y,z}(a_1,\dotsc,a_k)f^{\bigodot}_{x+1,y-1,z-1}(a_1,\dotsc,a_k)
\tag\edc
\endalign
$$
where the superscripts at the $f$'s emphasize the type of $F$-cored hexagon these formulas correspond to, and $x$, $y$ and $z$ have the same parity. 

By (\eda) and (\edb), what (\edc) states is that
$$
\spreadlines{3\jot}
\align
&
\M(FC^{\bigodot}_{x,y,z}(o,e))\M(FC^{\nwarrow}_{x,y-1,z-1}(o,e))\,
g_{x,y,z}(a_1,\dotsc,a_k)\,g_{x,y-1,z-1}(a_1,\dotsc,a_k)=
\\
&\ \ \ \ \ 
%=
\M(FC^{\nwarrow}_{x,y-1,z}(o,e))\M(FC^{\swarrow}_{x,y,z-1}(o,e))\,
g_{x,y-1,z}(a_1,\dotsc,a_k)\,g_{x,y,z-1}(a_1,\dotsc,a_k)
\\
&\!
+
\M(FC^{\leftarrow}_{x-1,y,z}(o,e))\M(FC^{\bigodot}_{x+1,y-1,z-1}(o,e))\,
g_{x-1,y,z}(a_1,\dotsc,a_k)\,g_{x+1,y-1,z-1}(a_1,\dotsc,a_k).
\\
\tag\edd
\endalign
$$

%!! NOTE: (3.2) and (3.3) is THE SAME recurrence!!
%
%   ALSO: (3.4) and (3.6) (resp., their close analogs (3.5) and (3.7)) are
% obtained from one another by swapping x and y !!
%
% SO we really only need to check THREE recurrences, and not six (and the last
% two should be very similar).

By (\edb), for fixed $a_1,\dotsc,a_k$, the value of $g_{x,y,z}(a_1,\dotsc,a_k)$ depends only on the sums $x+y$ and $x+z$. The values of these pairs for the six  occurrences of $g$ in (\edd) are
$$
\matrix
x+y & x+y-1 & | & x+y-1 & x+y & | & x+y-1 & x+y\\
x+z & x+z-1 & | & x+z & x+z-1 & | & x+z-1 & x+z,\\
\endmatrix
\tag\ede
$$
where we separated by vertical bars the quantities corresponding to the three terms of (\edd). In particular, for both rows, the values in the three portions into which they are broken by the vertical lines form the same set. By the observation at the beginning of this paragraph, it follows that the three products of $g$-functions in (\edd) have the same value.

Therefore, verifying (\edd) amounts to checking that
$$
\spreadlines{3\jot}
\align
&
\M(FC^{\bigodot}_{x,y,z}(o,e))\M(FC^{\nwarrow}_{x,y-1,z-1}(o,e))=
%\\
%&\ \ \ \ \ 
%=
\M(FC^{\nwarrow}_{x,y-1,z}(o,e))\M(FC^{\swarrow}_{x,y,z-1}(o,e))
\\
&\ \ \ \ \ \ \ \ \ \ \ \ \ \ \ \ \ \ \ \ \ \ \ \ 
+\M(FC^{\leftarrow}_{x-1,y,z}(o,e))\M(FC^{\bigodot}_{x+1,y-1,z-1}(o,e))
.
\tag\edf
\endalign
$$
However, this holds because it is a special case of the general recurrence\footnote{ This is so because $x$, $y$ and $z$ have the same parity; this is the only place in the verification of (\edc) where we use this assumption.} (\ecb)! This completes the proof of (\edc), and therefore that the $f_{x,y,z}(a_1,\dotsc,a_k)$'s satisfy recurrence~(\ecb).

In fact, as it turns out, this also proves that they satisfy recurrence (\ecc). Indeed, quite remarkably, upon deletion of the marks at superscripts (which recall are only used for geometric convenience, as each type is determined by the relative parities of the $x,y,z$ parameters), recurrences (\ecb) and~(\ecc) become identical. 

The same observation unifies (\ecd) and (\ecf), and also (\ece) and (\ecg). Thus to complete the proof it suffices to verify that the $f_{x,y,z}(a_1,\dotsc,a_k)$'s satisfy recurrences (\ecd) and (\ece). 

We consider first recurrence (\ecd). To show that the $f_{x,y,z}(a_1,\dotsc,a_k)$'s satisfy it amounts to verifying that
$$
\spreadlines{3\jot}
\align
&
f^{\swarrow}_{x,y,z}(a_1,\dotsc,a_k)f^{\bigodot}_{x-1,y-1,z}(a_1,\dotsc,a_k)=
%\\
%&\ \ \ \ \ \ \ \ \ \ \ \ \ \ \ \ \ \ \ \ \ \ \ \ \ \ \ 
%=
f^{\leftarrow}_{x,y-1,z}(a_1,\dotsc,a_k)f^{\nwarrow}_{x-1,y,z}(a_k,\dotsc,a_1)
\\
&\ \ \ \ \ \ \ \ \ \ \ \ \ \ \ \ \ \ \ \ \ \ \ \ \ \ \ \ \ \ \ \ \ \ \ \ \ \ \ 
+
f^{\bigodot}_{x,y,z-1}(a_1,\dotsc,a_k)f^{\swarrow}_{x-1,y-1,z+1}(a_1,\dotsc,a_k),
\tag\edg
\endalign
$$
for $x$ and $y$ of the same parity, $z$ of the opposite parity, and $k$ even. Using (\eda) and (\edb), (\edg) becomes
$$
\spreadlines{3\jot}
\align
&
\M(FC^{\swarrow}_{x,y,z}(o,e))\M(FC^{\bigodot}_{x-1,y-1,z}(o,e))\,
g_{x,y,z}(a_1,\dotsc,a_k)\,g_{x-1,y-1,z}(a_1,\dotsc,a_k)=
\\
&\ \ \ \ \ 
%=
\M(FC^{\leftarrow}_{x,y-1,z}(o,e))\M(FC^{\nwarrow}_{x-1,y,z}(e,o))\,
g_{x,y-1,z}(a_1,\dotsc,a_k)\,g_{x-1,y,z}(a_k,\dotsc,a_1)
\\
&\!
+
\M(FC^{\bigodot}_{x,y,z-1}(o,e))\M(FC^{\swarrow}_{x-1,y-1,z+1}(o,e))\,
g_{x,y,z-1}(a_1,\dotsc,a_k)\,g_{x-1,y-1,z+1}(a_1,\dotsc,a_k).
\\
\tag\edh
\endalign
$$
For the same reason as in the verification of (\ecb), the product of the two $g$-functions on the left hand side above is the same as the product of the last two $g$-functions on the right. 

For a different reason, it turns out that the product of the {\it first two} $g$-functions on the right in (\edh) is also equal to the product of the two $g$-functions on its left hand side! Indeed, even though the values of the $(x+y)$- and $(x+z)$-parameters of the two $g$-functions in the first term on the right hand side in (\edh) do not form the same set (namely, $\{x+y-2,x+y\}$, resp. $\{x+z-1,x+z\}$) as in the other two terms, it is not hard to verify using their explicit formula (\edb) that, for $x$ and $y$ of parity opposite to the parity of $z$, the $g$-functions satisfy the identity
$$
g_{x,y-1,z}(a_1,\dotsc,a_k)\,g_{x-1,y,z}(a_k,\dotsc,a_1)
=
g_{x,y,z}(a_1,\dotsc,a_k)\,g_{x-1,y-1,z}(a_1,\dotsc,a_k).
\tag\edi
$$
Therefore, dividing through equation (\edh) by the product of the two $g$-functions on the left, we obtain that (\edh) holds if and only if 
$$
\spreadlines{3\jot}
\align
&
\M(FC^{\swarrow}_{x,y,z}(o,e))\M(FC^{\bigodot}_{x-1,y-1,z}(o,e))
%\\
%&\ \ \ \ \ 
=
\M(FC^{\leftarrow}_{x,y-1,z}(o,e))\M(FC^{\nwarrow}_{x-1,y,z}(e,o))
\\
&\ \ \ \ \ \ \ \ \ \ \ \ \ \ \ \ \ \ \ \ \ \ \ \ \ \ \ \ \ \
+
\M(FC^{\bigodot}_{x,y,z-1}(o,e))\M(FC^{\swarrow}_{x-1,y-1,z+1}(o,e)).
\tag\edj
\endalign
$$
However, just as it was the case for the verification of recurrence (\ecb), this equality holds since it is a special case of the general recurrence (\ecd). We point out that here it is essential that the number of lobes $k$ is assumed to be even. Indeed, for $k$ odd, the 2-lobe $F$-cored hexagon entering the expression of $\M(FC^{\nwarrow}_{x-1,y,z}(a_k,\dotsc,a_1))$ via (\eda) would have, from left to right, lobes of sizes $o,e$, as opposed to $e,o$, as it is required by recurrence (\ecd). This completes veriying that the $f$'s satisfy recurrence (\ecd). 

Finally, we check that when the number of lobes $k$ is odd, the $f$'s satisfy recurrence~(\ece). For this, we need to verify that
$$
\spreadlines{3\jot}
\align
&
f^{\swarrow}_{x,y,z}(a_1,\dotsc,a_k)f^{\bigodot}_{x-1,y-1,z}(a_1,\dotsc,a_k)=
%\\
%&\ \ \ \ \ \ \ \ \ \ \ \ \ \ \ \ \ \ \ \ \ \ \ \ \ \ \ 
%=
f^{\leftarrow}_{x,y-1,z}(a_1,\dotsc,a_k)f^{\swarrow}_{x-1,y,z}(a_k,\dotsc,a_1)
\\
&\ \ \ \ \ \ \ \ \ \ \ \ \ \ \ \ \ \ \ \ \ \ \ \ \ \ \ \ \ \ \ \ \ \ \ \ \ \ \ 
+
f^{\bigodot}_{x,y,z-1}(a_1,\dotsc,a_k)f^{\swarrow}_{x-1,y-1,z+1}(a_1,\dotsc,a_k),
\tag\edk
\endalign
$$
for $x$ and $y$ of the same parity, $z$ of the opposite parity, and $k$ odd. Using (\eda) and (\edb), (\edk) becomes
$$
\spreadlines{3\jot}
\align
&
\M(FC^{\swarrow}_{x,y,z}(o,e))\M(FC^{\bigodot}_{x-1,y-1,z}(o,e))\,
g_{x,y,z}(a_1,\dotsc,a_k)\,g_{x-1,y-1,z}(a_1,\dotsc,a_k)=
\\
&\ \ \ \ \ 
%=
\M(FC^{\leftarrow}_{x,y-1,z}(o,e))\M(FC^{\swarrow}_{x-1,y,z}(o,e))\,
g_{x,y-1,z}(a_1,\dotsc,a_k)\,g_{x-1,y,z}(a_k,\dotsc,a_1)
\\
&\!
+
\M(FC^{\bigodot}_{x,y,z-1}(o,e))\M(FC^{\swarrow}_{x-1,y-1,z+1}(o,e))\,
g_{x,y,z-1}(a_1,\dotsc,a_k)\,g_{x-1,y-1,z+1}(a_1,\dotsc,a_k)
\\
\tag\edl
\endalign
$$
(note that unlike in (\edh), here the fourth 2-lobe $F$-cored hexagon has lobes listed in the order $(o,e)$, just like the other five; this is because for an odd number of lobes, the sum of the sizes of the odd indexed ones is the same if one starts from the beginning or from the end, and similarly for the even indexed lobes).

For the same reason as before, the product of the two $g$-functions on the left is the same as the product of the last two $g$-functions on the right. However, now the product of the first two $g$-functions on the right turns out to have a different value from these. Dividing through by the left hand side and using (\edb), (\ebga) and (\ebha), after simplifications (\edl) becomes
$$
1=\frac{x+y+z}{x+y+z+2o+2e-1}
+
\frac{2o+2e-1}{x+y+z+2o+2e-1},
\tag\edm
$$
which is indeed true. This completes the verification of the fact that the $f$-functions satisfy recurrence (\ece), and with this all recurrences (\ecb)--(\ecg) are checked.

\mysec{5. Proof of Theorem {\taa}. Geometrical interpretation}

{\it Proof of Theorem {\eaa}.} By the definition (\eac), we have
$$
\frac{\M(F^*(a_1,\dotsc,a_k))}{\M(F^*(o,e))}
=
\lim_{N\to\infty}\frac{\M(FC_{N,N,N}(a_1,\dotsc,a_k))}{\M(FC_{N,N,N}(o,e))}
\tag\eea
$$
The fraction whose limit is taken above is expressed by formula (\ebb), which dramatically simplifies --- to just $s(a_1,\dotsc,a_{k-1})s(a_2,\dotsc,a_k)$ --- whenever the $y$- and $z$-parameters have the same value. Since this is the case in (\eea) (both have value $N$), we obtain the statement of Theorem {\eaa}. \epf
%Remark --- enough $x=y$, get same result! So it's same for a family of diifferent Gibbs measures. The efects in numerator and denominator cancel out precisely.

%Geom. int.

\topinsert
\centerline{\mypic{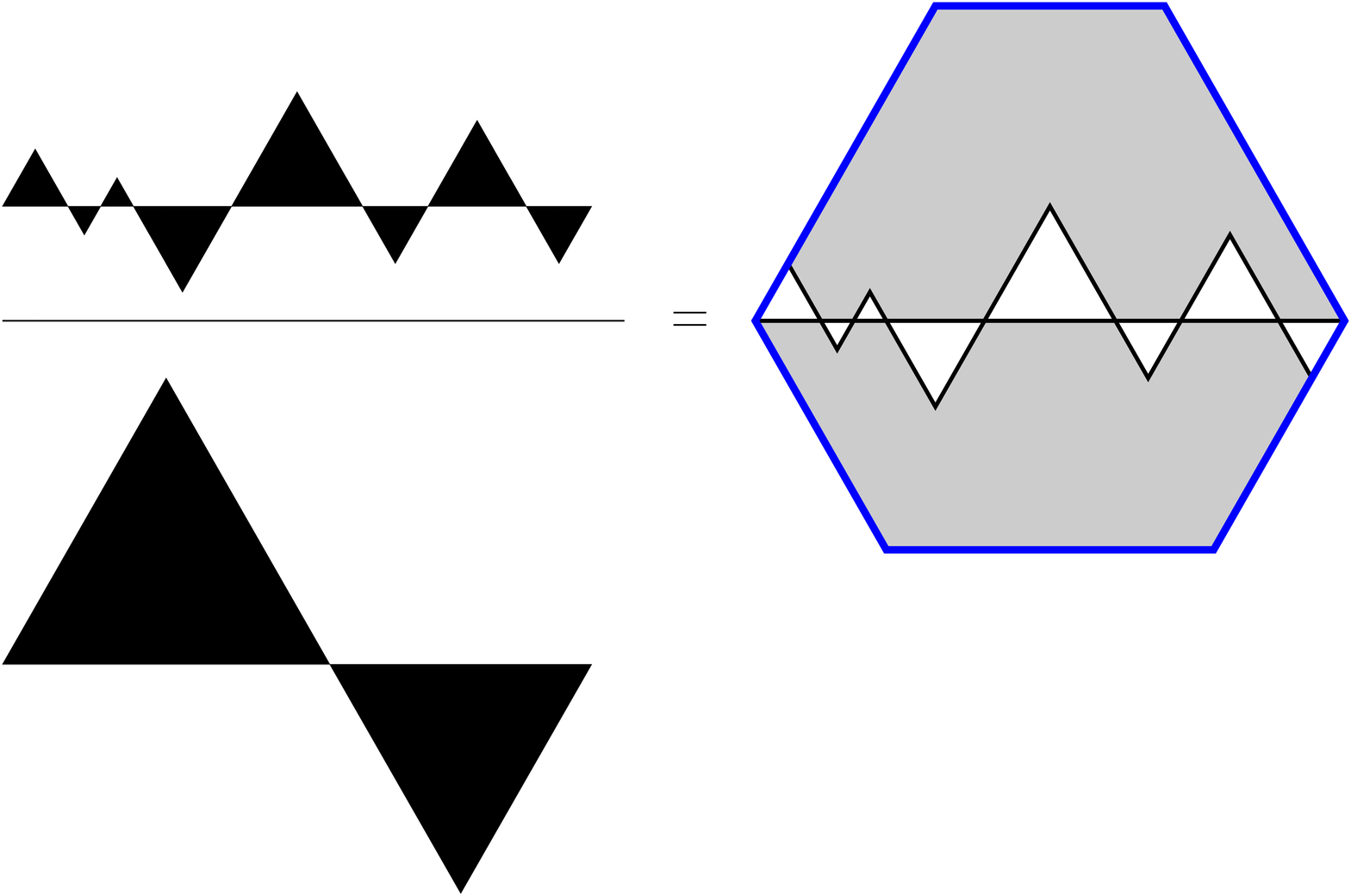}}
\medskip
\centerline{{\smc Figure~{\fea}. {\rm Geometrical interpretation of Theorem {\taa}.}}}
\endinsert

{\bf Remark 4.}
It turns out that Theorem {\taa} has an especially simple geometrical interpretation: 
$$
\frac{\M(F^*(a_1,\dotsc,a_k))}{\M(F^*(o,e))}
=
\M(H_F),
\tag\eeb
$$
where $H_F$ is the smallest hexagon on the triangular lattice that contains the fern and has the property that the region obtained from it by removing the fern is balanced (i.e., contains the same number of up-pointing and down-pointing unit triangles). This is illustrated in Figure {\fea}. 

Indeed, directly from Theorem {\taa} we have that the left hand side of (\eeb) is equal to the product of the numbers of tilings of the semihexagons $S(a_1,\dotsc,a_{k-1})$ and $S(a_2,\dotsc,a_{k})$ (see Figure {\fad} for the definition of such a semihexagon). However, due to forced tiles, these semihexagons can be augmented to the regions above and below the fern in $H_F$, respectively, without affecting the number of their tilings (when the number of lobes is even, as in Figure {\fea}, each semihexagon gets augmented by a rhombus of forced tiles; when the number of lobes is odd, only the bottom semihexagon gets augmented, by two rhombi of forced lozenges). Since the latter two regions are necessarily tiled internally in any tiling of $H_F$, $\M(H_F)$ is equal to the product of their numbers of tilings, and we obtain~(\eeb).

%This raises a compelling question: Is perhaps analogous ratio for shamrocks also equal to M(enveloping hexagon)? There is a chance that it is so --- in \cite{\vf} we normalized by  F(a+b+c,m), which is what we need in this context. So question is:
%                   ?
%P(a,b,m)P(a+b,c,m) = P(a,b,m)P(b,c,m)P(a,c,m) ?

%And the answer is no! So for shamrocks it is not true. 

\medskip
{\bf Remark 5.} We note that in fact, by formula (\ebb), the value of the limit on the right hand side of (\eea) stays the same even if the list of indices in the $F$-cored hexagons are $x,y,y$ instead of $N,N,N$ --- for {\it arbitrary} values of $x$ and $y$ (even not necessarily approaching infinity)! 

In particular this is one way to see (\eeb) directly: simply use that, by the previous observation, the right hand side of (\eea) is equal to its $x=y=z=0$ specialization. Then~(\eeb) follows by noting that for $x=y=z=0$, the numerator of the fraction on the right hand side of (\eea) becomes $\M(H_F)$, while the denominator becomes the number of lozenge tilings of the smallest lattice hexagon containing $F(o,e)$, which is 1.

\mysec{6. Concluding remarks} 
%and an open problem}

In this paper we introduced a counterpart structure to the shamrocks studied in \cite{\vf}, which, just like the latter, can be included at the center of a lattice hexagon on the triangular lattice so that the region obtained from the hexagon by removing it has its number of lozenge tilings given by a simple product formula. The new structure, called a fern, consists of an arbitrary number of equilateral triangles of alternating orientations lined up along a lattice line. The shamrock and the fern seem to be the only structures with this property. It would be interesting to understand why these are the only two such structures.

\mysec{References}
{\openup 1\jot \frenchspacing\raggedbottom
\roster

\myref{\And}
  G. E. Andrews, Plane partitions (III): The weak Macdonald
conjecture, {\it Invent. Math.} {\bf 53} (1979), 193--225.

\myref{\cekz}
  M. Ciucu, Eisenkolbl, C. Krattenthaler and D. Zare,
Enumeration of lozenge tilings of hexagons with a central triangular hole,
{\it J. Combin. Theory Ser. A} {\bf 95} (2001), 251--334.

%\myref{\sc}
%  M. Ciucu, A random tiling model for two dimensional electrostatics, 
%{\it Mem. Amer. Math. Soc.} {\bf 178} (2005), no. 839, 1--106.

\myref{\ec}
  M. Ciucu, The scaling limit of the correlation of holes on the triangular lattice
with periodic boundary conditions, {\it Mem. Amer. Math. Soc.} {\bf 199} (2009),
no. 935, 1-100.

\myref{\ov}
  M. Ciucu, Dimer packings with gaps and electrostatics, {\it Proc. Natl. Acad. Sci.
USA} {\bf 105} (2008), 2766-2772.

%\myref{\ef}
%  M. Ciucu, The emergence of the electrostatic field as a Feynman sum in random
%tilings with holes, {\it Trans. Amer. Math. Soc.} {\bf 362} (2010), 4921-4954.

\myref{\vf}
  M. Ciucu and C. Krattenthaler,
A dual of MacMahon's theorem on plane partitions,
{\it Proc. Natl. Acad. Sci. USA} {\bf 110} (2013), 4518--4523.

\myref{\gd}
  M. Ciucu, The interaction of collinear gaps of arbitrary charge in a two
dimensional dimer system, {\it Comm. Math. Phys.}, {\bf 330} (2014), 1115--1153.

%\myref{\anglepap}
%  M. Ciucu and I. Fischer, A triangular gap of size two in a sea of dimers on a
%$60^\circ$ angle, preprint, July 2012, submitted.

\myref{\CLP}
  H. Cohn, M. Larsen, and J. Propp, The shape of a typical boxed plane partition,
{\it New York J. of Math.} {\bf 4} (1998), 137--165.

%\myref{\DT}
%  G. David and C. Tomei, The problem of the calissons,
%{\it Amer\. Math\. Monthly} {\bf 96} (1989), 429--431.

\myref{\GT}
  I. M. Gelfand and M. L. Tsetlin, Finite-dimensional representations of the group of unimodular matrices (in Russian), {\it Doklady Akad. Nauk. SSSR (N. S.)} {\bf 71} (1950), 825--828.

%\myref{\Glaish}
%  J. W. L. Glaisher, On Certain Numerical Products in which the Exponents Depend Upon the Numbers, {\it  Messenger Math.} {\bf 23} (1893), 145--175.

\myref{\KKZ}
  C. Koutschan, M. Kauers and D. Zeilberger, A proof of George
Andrews' and David Robbins' $q$-TSPP-conjecture, {\it Proc Natl Acad Sci USA} {\bf 108} (2011), 2196--2199.

\myref{\Kuo}
E. H. Kuo, Applications of graphical condensation for enumerating matchings and tilings. Theoret. Comput. Sci. 319 (2004), no. 1-3, 29–57.

\myref{\Kup}
  G. Kuperberg, Symmetries of plane partitions and the permanent-de\-ter\-mi\-nant
method, {\it J. Combin. Theory Ser. A} {\bf 68} (1994), 115--151.

\myref{\MacM}
  P. A. MacMahon, Memoir on the theory of the partition of numbers---Part V. 
Partitions in two-dimensional space, {\it Phil. Trans. R. S.}, 1911, A.

\myref{\Sta}
  R. P. Stanley, Symmetries of plane partitions,
{\it J. Comb. Theory Ser. A} {\bf 43} (1986), 103--113.

\myref{\Ste}
  J. R. Stembridge, The enumeration of totally symmetric plane partitions,
{\it Adv. in Math.}, {\bf 111} (1995), 227--243.

\endroster\par}

\enddocument